\documentclass[11pt]{article}
\usepackage{graphicx}
\usepackage{amssymb}
\usepackage{epstopdf}
\usepackage{amsmath}
\usepackage{amsfonts}
\newcommand{\braket}[2]{\langle #1,#2 \rangle}
\newcommand{\la}{\lambda}
\newcommand{\var}{\varepsilon}

\DeclareSymbolFont{AMSb}{U}{msb}{m}{n}
\DeclareMathSymbol{\N}{\mathbin}{AMSb}{"4E}
\DeclareMathSymbol{\Z}{\mathbin}{AMSb}{"5A}
\DeclareMathSymbol{\R}{\mathbin}{AMSb}{"52}
\DeclareMathSymbol{\Q}{\mathbin}{AMSb}{"51}
\DeclareMathSymbol{\I}{\mathbin}{AMSb}{"49}
\DeclareMathSymbol{\C}{\mathbin}{AMSb}{"43}

\begin{document}
 
\addtolength{\textheight}{0 cm}
\addtolength{\hoffset}{0 cm}
\addtolength{\textwidth}{0 cm}
\addtolength{\voffset}{0 cm}

\setcounter{secnumdepth}{5}
 \newtheorem{proposition}{Proposition}[section]
\newtheorem{theorem}{Theorem}[section]
\newtheorem{lemma}[theorem]{Lemma}
\newtheorem{coro}[theorem]{Corollary}
\newtheorem{remark}[theorem]{Remark}
\newtheorem{ex}[theorem]{Example}
\newtheorem{claim}[theorem]{Claim}
\newtheorem{conj}[theorem]{Conjecture}
\newtheorem{definition}[theorem]{Definition}
\newtheorem{application}{Application}
 
\newtheorem{corollary}[theorem]{Corollary}
\def\HADX{{\cal H}_{\rm AD}(X)}
\def\HADY{{\cal H}_{\rm AD}(Y)}
\def\HADH{{\cal H}_{\rm AD}(H)}
\def\HTADX{{\cal H}_{\rm TAD}(X)}
\def\HTADY{{\cal H}_{\rm TAD}(Y)}
\def\HTADH{{\cal H}_{\rm TAD}(H)}
\def\LX{{\cal L}(X)}
\def\LY{{\cal L}(Y)}
\def\LH{{\cal L}(H)}
 \def\ASD{{\cal L}_{\rm AD}(X)}
 \def\ASDY{{\cal L}_{\rm AD}(Y)}
\def\ASDH{{\cal L}_{\rm AD}(H)}
 \def\ASDP{{\cal L}^{+}_{\rm AD}(X)}
  \def\ASDYP{{\cal L}^{+}_{\rm AD}(Y)}
   \def\ASDHP{{\cal L}^{+}_{\rm AD}(H)}
    \def\TADX{{\cal L}_{\rm TAD}(X)}
        \def\TADY{{\cal L}_{\rm TAD}(Y)}
            \def\TADH{{\cal L}_{\rm TAD}(H)}
 \def\CX{{\cal C}(X)}
\def\CY{{\cal C}(Y)}
\def\CH{{\cal C}(H)}
 
\def\PX{{\cal A}(X)}
\def\PY{{\cal A}(Y)}
\def\PH{{\cal A}(H)}
\def\phi{{\varphi}}
\def\AH{A^{2}_{H}}

\title{Anti-selfdual Hamiltonians: Variational resolutions for Navier-Stokes and other nonlinear  evolutions}
\author{ Nassif  Ghoussoub\thanks{Research partially supported by a grant
from the Natural Sciences and Engineering Research Council of Canada. The
author gratefully acknowledges the hospitality and support of the Centre de Recherches Mathematiques and the Universit\'e de Nice where this work  was initiated.} 
\\
\small Department of Mathematics,
\small University of British Columbia, \\
\small Vancouver BC Canada V6T 1Z2 \\
\small {\tt nassif@math.ubc.ca} 
\\
}
\maketitle

\section*{Abstract}  The theory of anti-selfdual (ASD) Lagrangians developed  in \cite{G2} allows a variational resolution for equations of the form $\Lambda u+Au +\partial \phi (u)+f=0$ where $\phi$ is a convex lower-semi-continuous function on a  reflexive Banach space $X$, $f\in X^*$, $A: D(A)\subset X\to X^*$ is a positive  linear operator and where $\Lambda: D(\Lambda)\subset X\to X^{*}$ is a non-linear operator that satisfies suitable continuity and anti-symmetry properties. ASD Lagrangians on path spaces also yield variational resolutions for nonlinear evolution equations of the form 
$\dot u (t)+\Lambda u(t)+Au(t) +f\in -\partial \phi (u(t))$ starting at $u(0)=u_{0}$. In both  stationary and  dynamic cases, the equations associated to  the proposed variational principles  are not derived from the fact they are critical points of the action functional, but because they are also zeroes of the Lagrangian  itself.The approach has many applications, in particular to Navier-Stokes type equations and to the differential systems of hydrodynamics, magnetohydrodynamics and thermohydraulics. 
 
\tableofcontents

\newpage

\section{Introduction} A new variational framework was developed in \cite{G2} where solutions of various equations, not normally of Euler-Lagrange type, can still be obtained as minima of functionals of the form
\[
 I(u)=L(u, A u)+\ell (b_1(x), b_2(x)) \quad  \hbox{\rm or \quad $I(u)=\int_{0}^{T}L(t, u(t), \dot u(t)+A u(t))dt+\ell (u(0), u(T)).$}
 \]
The Lagrangians $L$ (and $\ell$) must obey certain anti-selfdual conditions, while the operators $A$ are  essentially  skew-adjoint  modulo boundary terms represented by a pair of  operators $(b_1,b_2)$. For such ``anti-selfdual" (ASD) Lagrangians,  the minimal value will always be zero and --just like the self (and antiself) dual equations of  quantum field theory (e.g. Yang-Mills and others)-- the equations associated to such minima are not derived from the fact they are critical points of the functional $I$, but because they are also zeroes of the Lagrangian $L$ itself. In other words, the solutions will satisfy
  \[
  L(u, A u)+\langle u,A u\rangle=0 \quad {\rm and} \quad L(t, u(t), \dot u(t)+A_{t}u(t))+ \langle u (t),\dot u (t) \rangle=0.
  \]
 It is also shown in \cite{G2} that ASD Lagrangians possess remarkable  permanence properties making them more  prevalent than expected and quite easy to construct and/or identify. The variational game changes from  the  analytical proofs of existence of extremals for general action functionals, to a more  algebraic search of an appropriate ASD Lagrangian for which the minimization problem is remarkably simple with value always equal to zero. This makes  them efficient new tools for proving existence and  uniqueness results for a large array of differential equations.  

In this paper, we tackle  boundary value problems of the form:
 \begin{equation}
 \left\{ \begin{array}{lcl}
\label{eqn:laxmilgram1000}
 \hfill - \Lambda u-Au+f&\in &\partial \varphi (u)\\
\hfill  b_1(u)&=&0 \\
\end{array}\right.
\end{equation}
  as well as  parabolic evolution equations of the form:
 \begin{equation}
 \left\{ \begin{array}{lcl}
\label{eqn:1001}
 \hfill   -\dot{u}(t)  -\Lambda u (t)-A u(t)&\in &\partial \varphi (t, u(t)) \quad \quad \hbox {\rm a.e. $t\in [0, T]$}\\
\hfill b_{1}(u(t))&=&b_1(u_{0}) \quad \quad \quad \quad \hbox{\rm a.e $t\in [0, T]$}\\
\hfill u(0)&=&u_{0}
\end{array}\right.
\end{equation}
where $u_{0}$ is a given initial value. Here $\phi$ is a convex lower semicontinuous functional, 
 $\Lambda$ is a  non-linear ``conservative" operator,  $A$ is a linear not necessarily bounded but essentially skew-adjoint operator modulo  the operators $(b_1,b_2)$. 
  
 As applications to our method, we provide a variational resolution to equations involving  nonlinear operators such as the Navier-Stokes equation for a fluid driven by its boundary:
  \begin{equation}
\label{NSE3}
 \left\{ \begin{array}{lcl}
    \hfill
 (u\cdot \nabla)u +f &=&\nu \Delta u - \nabla  p \quad \hbox{\rm on $ \Omega$}\\
\hfill {\rm div} u&=&0 \quad \hbox{\rm on  $\Omega$}\\
\hfill u&=&u^0 \quad \hbox{\rm on $\partial \Omega$}\\
\end{array}\right.\nonumber 
\end{equation}
where $u^{0} \in H^{3/2}(\partial \Omega)$ is such that $\int_{\partial \Omega} u^{0} {\bf \cdot n}\, d\sigma =0$, $\nu>0$ and $f\in L^{p}(\Omega;\R^{3})$. \\
We can also deal with the superposition of such non-linear operators with 
 non self-adjoint first order operators such as  linear transport maps:
\begin{equation}
\label{NSE5}
 \left\{ \begin{array}{lcl}
    \hfill
 (u\cdot \nabla)u + \vec a\cdot  \nabla  u + a_0 u+|u|^{m-2}u+f &=&\nu \Delta u - \nabla  p \quad \hbox{\rm on $ \Omega$}\\
\hfill {\rm div} u&=&0 \quad \hbox{\rm on  $\Omega$}\\
\hfill u &=&0 \quad \hbox{\rm on $\partial \Omega$}\\
\end{array}\right.\nonumber 
\end{equation}
where $\vec a \in C^{\infty}(\bar\Omega)$ is a smooth vector field and  $a_0 \in L^{\infty}$  are such that  $a_0-\frac{1}{2} {\rm div}  (a)  \geq 0$.\\
The methods extend to the dynamic case where typically we give a variational resolution to the Navier-Stokes evolution
   \begin{equation}
\label{NSE}
 \left\{ \begin{array}{lcl}
    \hfill
\frac{\partial u}{\partial t}+(u\cdot \nabla)u -f &=&\nu \Delta u - \nabla  p \quad \hbox{\rm on $[0,T]\times \Omega$}\\
\hfill {\rm div} u&=&0 \quad \hbox{\rm on $[0,T]\times \Omega$}\\
\hfill u(t,x)&=&0 \quad \hbox{\rm on $[0,T]\times \partial \Omega$}\\
\hfill u(0,x)&=&u_{0}(x) \quad \hbox{\rm on $\Omega$.} 
\end{array}\right.\nonumber
\end{equation}

 The paper, though sufficiently self-contained, is better read in conjunction with \cite{G2}. It is  organized as follows: In section 2, we introduce the concept of anti-selfdual Hamiltonian which is the appropriate dual notion to anti-selfdual Lagrangians.  In section 3, we give the main non-linear variational principle, which is applied in section 4  to obtain variational proofs for the existence of  stationary solutions for various nonlinear equations of Lax-Milgram type. In section 5,  we deal with the dynamic case where we  provide a variational resolution to several nonlinear parabolic initial-value problems, including those  appearing in the basic models hydro-dynamics. 
 
\section{Basic properties of Anti-selfdual Hamiltonians}

\begin{definition} Let $X$ be a reflexive Banach space.  Say that  a functional $H:X\times X\to \R\cup\{+\infty\}\cup\{-\infty\}$ is an {\it anti-selfdual Hamiltonian} if for each $y\in X$, the function $x\to -H(x,y)$ from $X$ to $\R\cup\{+\infty\}\cup\{-\infty\}$ is convex and the function $x\to H(-y,-x)$ is its convex lower semi-continuous envelope.
\end{definition}
It readily follows that for an ASD Hamiltonian $H$, the function  $y\to H(x,y)$ is convex and lower semi-continuous for each $x\in X$, and that  the following inequality holds for every $(x,y)\in X\times X$,
 \begin{equation}
\label{almost.odd}
H(-y,-x)\leq -H(x,y),
\end{equation}
In particular, we have for every $x\in X$,
 \begin{equation}
\label{negative.diagonal}
H(x,-x)\leq 0.
\end{equation}
The class of  anti-selfdual Hamiltonian on a space $X$, will be denoted by $\HADX$. The most basic ASD Hamiltonian is $H(x,y)=\|y\|^{2}-\|x\|^{2}$ (Maxwell's Hamiltonian)  or more generally $H(x,y)=\phi (-y)-\phi (x)$  where $\phi$ is any finite convex lower semi-continuous function on $X$. More generally,  if $B: X\to X^*$ is a skew-adjoint bounded linear operator, $f\in X^*$, and if $\phi:X\to \R\cup\{+\infty\}$ is proper convex and lower semi-continuous, then 
 \begin{equation}
  H(x,y)= \left\{ \begin{array}{lcl}
  \phi (-y)-\phi (x)-\langle  B x, y\rangle +\langle f, x+y\rangle &&\hbox{ if $x\in {\rm Dom}(\phi)$}\\
 -\infty &&\, \, \hbox{if $x\notin {\rm Dom} (\phi)$}
       \end{array}  \right.
    \end{equation}
  is also an anti-selfdual Hamiltonian. 
  We define the (partial) domain of $H$ to be
  \begin{equation}
{\rm Dom}_1(H)=\{x\in X; H(x,y)>-\infty \hbox{ for all $y\in X$} \}. 
  \end{equation}
  Note that if $\phi$ is a convex lower semi-continuous function that is bounded below on $X$, and if $H(x,y)=\phi (-y) -\phi (x)$ is the anti-selfdual  Hamiltonian  associated to $\phi$, then  ${\rm Dom}_1 H={\rm Dom}{\phi}$. 
Note also that for any $z\in  {\rm Dom}_1(H)$, we have that the function $\phi_z:x\to -H(x,-z)$ is convex and valued in  $\R\cup\{+\infty\}$. Moreover ${\rm Dom}_1(H) \subset {\rm Dom}(\phi_z)$, hence for any $z, y \in {\rm Dom}_1(H)$,
\[
\hbox{$H(z,-y)= -H(y,-z)$ if and only if the function $x\to H(x,-z)$ is upper semi-continuous at $y$.}
\]
We can now introduce the following
\begin{definition} Say that  an Anti-selfdual Hamiltonian $H:X\times X\to \R$ is {\it tempered} if for every $y\in {\rm Dom}_1(H)$,  the function $x\to H(x,-y)$ is concave and upper semi-continuous from $X$ to $\R\cup\{-\infty\}$.
\end{definition}
It then follows that 
\begin{equation}
\hbox{ $H(y,x)=-H(-x,-y)$  for all $(x,y)\in X\times {\rm Dom}_1(H)$.}
 \end{equation}
and therefore
\begin{equation}
\label{null.diagonal}
H(x,-x)=0 \quad \hbox{\rm for all $x\in {\rm Dom}_1(H)$}.
\end{equation}
The class of tempered anti-selfdual Hamiltonian on a space $X$, will be denoted by $\HTADX$.

The most basic tempered ASD Hamiltonian is 
 $
  H(x,y)=\phi (y)-\phi (-x)+\langle x,  B y\rangle +\langle f, x+y\rangle
   $
 where $\phi$ is any finite convex lower semi-continuous function on $X$,  $f\in X^*$, and where $ B:X\to X^*$ is a skew-adjoint bounded linear operator.
   Tempered ASD Hamiltonians satisfy some obvious  permanence properties that we summarize in the following proposition.
  \begin{proposition}
 Let $X$ be a reflexive Banach space, then the following holds:
   \begin{enumerate}
\item  If $H$ and $K$ are in $\HTADX$ and $\lambda >0$,  then the Hamiltonians $H+ K$ (defined as $-\infty$ if the first variable is not in  ${\rm Dom}_1(H)\cap  {\rm Dom}_1(K)$),  and $\lambda {\bf \cdot} H$  also belong to $\HTADX$. 
\item If $H_{i}\in {\cal H}_{\rm TAD}(X_{i})$ where $X_{i}$ is a reflexive Banach space for each $i\in I$, then the Hamiltonian $H:=\Sigma_{i\in I} H_{i}$ defined by 
$H((x_{i})_{i}, (y_{i})_{i})=\Sigma_{i\in I} H_{i}(x_{i}, y_{i})$
  is  in ${\cal H}_{\rm TAD}(\Pi_{i\in I} X_{i})$.
   \item If $H\in \HTADX$ and $ B :X\to X^{*}$ is a skew-adjoint bounded linear operator   then the Hamiltonian $H_{ B}$ defined by 
   $
   H_{ B}(x,y)=H(x,y)+\langle  B x, y\rangle$
    is also in $\HTADX$.
  \item If $H \in \HTADX$ and $K\in \HTADY$, then for any bounded linear operator $A: X\to Y^{*}$,  the Hamiltonian $H+_{A}K$ defined by 
\[
 (H+_{A}K)((x,y), (z,w))=H(x,z)+K(y, w)+\langle A^{*}y, z \rangle -\langle Ax, w \rangle
 \]
   belongs to ${\cal H}_{\rm TAD}(X\times Y)$.
 \item  If $\phi$ is a proper convex lower semi-continuous function on $X\times Y$ and $A$ is any bounded linear operator $A: X\to Y^{*}$,  then the Hamiltonian $H_{\phi, A}$ defined by 
 \[
 H_{\phi, A}((x,y), (z,w))=\phi (-z,-w)-\phi (x,y) +\langle A^{*}y, z \rangle -\langle Ax, w \rangle
 \]
 also  belongs to ${\cal H}_{\rm TAD}(X\times Y)$
 \end{enumerate}
\end{proposition}

   This notion is in a certain sense dual to the notion of anti-selfdual Lagrangian introduced and developed in \cite{G2}. Indeed, let  $\LX$ be the class of convex Lagrangians on a reflexive Banach space $X$. These are all functions  $L:X\times X^{*} \to \R \cup \{+\infty\}$ which are convex and lower semi-continuous (in both variables) and which are not identically $+ \infty$. The (partial) domain of the Lagrangian $L$ is defined as
  \begin{equation}
{\rm Dom}_1(L)=\{x\in X; L(x,p)<+\infty \hbox{ for some $p\in X^*$} \}. 
  \end{equation}  
  To each Lagrangian $L$ on $X\times X^{*}$, we can associate its {\it Hamiltonian} on $X\times X$ defined as the Legendre transform in the second variable, i.e.,
 \[
H_{L}(x,y):=\sup\{\langle p, y\rangle -L(x,p); p\in X^{*}\}.
\]
It is clear that   ${\rm Dom}_1(L)={\rm Dom}_1(H_L)$. \\
The Legendre-Fenchel dual (in both variables) of $L$ is defined at any pair $(q,y)\in X^{*}\times X$ by: 
 \[
    L^*( q,y)= \sup \{ \braket{q}{x} + \braket{y}{p} - L(x,p);\, x \in X, p \in X^{*}  \}.
\]
We recall the notion of anti-selfdual Lagrangians developed in \cite{G2}.
 \begin{definition} \rm Let $L$ be a Lagrangian in  $\LX$. We say that  \\
 (1) $L$ is  an {\it anti-self dual Lagrangian} on $X\times X^{*}$, if \\ 
\begin{equation}
\label{seldual}
L^*( p, x) =L(-x, -p ) \quad \hbox{\rm for all $(p,x)\in X^{*}\times X$}. 
\end{equation}
(2) $L$ is {\it anti-self dual on the graph of $\Lambda$}, the latter being a map from a subset $D\subset X$ into $X^{*}$, if 
\begin{equation}
\label{seldual.atgraph}
L^*(\Lambda x, x) =L(-x, -\Lambda x ) \quad \hbox{\rm for all $x\in D$}. 
\end{equation} 
\end{definition}
We denote by $\ASD$ the class of ASD-Lagrangians. We now procced to identify the class of Hamiltonians associated to  $ASD$-Lagrangians.  We denote by  $K_2^*$ (resp., $K_1^*$)  the Legendre dual of a functional $K(x,y)$ with respect to the second variable (resp., the first variable), we have the following
 
 \begin{proposition} Let $L$ be an ASD Lagrangian on a reflexive Banach space $X$, then its corresponding Hamiltonian $H=H_{L}$ is  anti-selfdual. 
\end{proposition}
\noindent {\bf Proof:} Since a Lagrangian $L\in \LX$ is convex in both variables, its corresponding Hamiltoninan $H_L$ is always concave in the first variable. Also note that  the Legendre transform of  $-H_L(\cdot, y)$ with respect to the first variable is related to the Legendre transform in both variables of its Lagrangian in the following way. 
\begin{eqnarray*}
(-H_L)^*_1(p,y)&=&\sup\{\langle p, x\rangle +H_L(x,y); x\in X\}\\
&=&\sup\{\langle p, x\rangle +\sup\{\langle y, q\rangle-L(x,q); q\in X^*\}; x\in X\}\\
&=&L^*(p,y).
\end{eqnarray*} 
If now $L$ is an ASD Lagrangian, then the convex lower semi-continuous envelope of the function $x\to -H_L(x, y)$ (i.e., the largest convex lower semi-continuous function below the function $x\to -H_L(x, y)$)  is
\begin{eqnarray*}
(-H_L)^{**}_1(x,y)&=&\sup\{\langle p, x\rangle -(-H_L)^*_1(p,y); p\in X^*\}\\
&=&\sup\{\langle p, x\rangle  -L^*(p,y); p\in X^*\}\\
&=&\sup\{\langle p, x\rangle  -L(-y,-p); p\in X^*\}\\
&=&H_L(-y,-x).
\end{eqnarray*} 
Note that  a characterization of anti-selfdual Hamiltonian that correspond to an ASD Lagrangian (i.e., $H=H_{L}$ for some $L\in \ASD$) is that 
\[
 H_{2}^{*}(-x,-p)=(-H)_{1}^{*}(p,x),
 \]
 for each $(x,p)\in X\times X^{*}$. In this case,  the corresponding ASD Lagrangian is nothing else but $L(x,p):=H_{2}^{*}(x,p)=(H)_{1}^{*}(-p,-x)$. \\

As mentioned above  since a Lagrangian $L\in \LX$ is convex in both variables, then its corresponding Hamiltoninan $H_L$ is always concave in the first variable. However,  $H_L$ is not necessarily upper semi-continuous in the first variable, even if $L$ is an anti-selfdual Lagrangian. This leads to the following notion. 
 
\begin{definition} A Lagrangian $L\in \LX$ will be called {\it tempered} 
if for each  $y\in {\rm Dom}_1(H)$, the map $x\to H(x,-y)$ from  $X$ to $\R \cup \{-\infty\}$ is upper semi-continuous.  
 \end{definition}
 A typical tempered Lagrangian (resp., tempered ASD-Lagrangian) is $L(x,p)=\phi (x) +\psi^{*}(p)$ (resp., $L(x,p)=\phi (x) +\phi^{*}(-p)$) where $\phi$ and $\psi$ are convex and lower semi-continuous on $X$.  We let ${\cal L}_{_{T}}(X)$ denote the class of tempered Lagrangians and ${\cal L}_{_{TAD}}(X)$ the class of tempered ASD Lagrangians on $X$.

 We now recall from \cite{G2} a few of the operations defined on the class of Lagrangians $\LX$ and study the permanence properties of the class $\TADX$  of tempered ASD Lagrangians. 
  \begin{itemize}
 \item {\it Addition:} If $L, M\in \LX$, define the Lagrangian $L\oplus M$ on $X\times X^{*}$ by:
 \[
 (L\oplus M)(x,p)=\inf\{L(x, r) + M(x,p-r); r\in X^{*}\}
 \]
\item {\it Convolution:} If $L, M\in \LX$, define the Lagrangian $L\star M$ on $X\times X^{*}$ by:
\[
(L\star M) (x,p)=\inf\{L(z, p) + M(x-z,p); z\in X\}
\]
\item {\it Right operator shift:}   If $L\in \LX$ and $ B :X\to X^{*}$ is a bounded linear operator, define the Lagrangian $L_{ B}$ on $X\times X^{*}$ by
\[
L_{ B} (x,p):=L(x,  B x+p).
\] 
  \end{itemize}
  \begin{lemma}  Let $X$ be a reflexive Banach space, then the following hold:
 \begin{enumerate}
 \item If  $L$ and $M$ are two Lagrangians in $\LX$, then $
 H_{L\oplus M}(x,y)=H_{L}(x,y)+H_{M}(x,y),
$ where $H_L$ and $H_M$ denote the corresponding Hamiltonians.
 \item If $L$ and $M$ are in $ \ASD$, then $L^{*}\oplus M^{*}(q,y)=L\star M (-y,-q)$ for every $(y,q)\in X\times X^{*}$. 
 \item If $L$ is an ASD Lagrangian and $M$ is of the form $M(x,p)=\phi (x)+\phi^*(-p)$ for some convex l.s.c. function $\phi$,   then $(L\oplus M)^{*} = L^{*}\star M^{*}$ and $(L\star M)^{*} = L^{*}\oplus M^{*}$.

 \end{enumerate}
 \end{lemma}
 \noindent{\bf Proof:}  
 (1) and (2) are straightforward, while (3) was established in \cite{G2}. It follows that the $\lambda$-regularization of an ASD Lagrangian $L$, that is 
$L_\lambda:=L\star T_\lambda$ where $T_\lambda(x,p)=\frac{\lambda^{2}\|x\|^{2}}{2}+\frac{\|p\|^{2}}{2\lambda^{2}}$, is also an ASD Lagrangian. \\

We shall see later that not all ASD Lagrangians are automatically tempered. This lemma shows that it is the case under certain coercivity conditions.
   \begin{proposition} Let $L$ be an ASD Lagrangian on a reflexive Banach space $X$. 
     If  for some $p_0\in X$ and  $\alpha >1$, we have that $L(x,p_0) \leq C(1+\|x\|^\alpha)$ for all $x\in X$, then $L$ belongs to  ${\cal L}_{_{TAD}}(X)$. 
   \end{proposition}
\noindent{\bf Proof:} Note that in this case,  we readily have that ${\rm Dom}_1(L)={\rm Dom}_1(H_L)=X. $\\
Assume first that $\lim\limits_{\|x\|+\|p\|\to +\infty}\frac{L(x,p)}{\|x\|+\|p\|}= \infty$, and write 
\begin{eqnarray*}
H_L(x,y)&=&\sup\{\langle p, y\rangle -L(x,p); p\in X^{*}\}\\
&=&\sup\{\langle p, y\rangle -L^{*}(-p,-x); p\in X^{*}\}\\
&=&\sup\{\langle p, y\rangle -\sup\{\langle -p, z\rangle +\langle -x, q\rangle- L(z,q); \, z\in X, q\in X^{*}\}; \, p\in X^{*}\}\\
&=&\sup\{\langle p, y\rangle +\inf\{\langle p, z\rangle +\langle x, q\rangle+ L(z,q); \, z\in X, q\in X^{*}\}; \, p\in X^{*}\}\\
&=&\sup\limits_{p\in X^{*}}\inf\limits _{(z,q)\in X\times X^{*}}\{\langle p, y\rangle +\langle p, z\rangle +\langle x, q\rangle+ L(z,q)\}.
\end{eqnarray*}
The function $S$ defined on the product space $(X\times X^{*})\times X^{*})$ as
\[
S((z,q), p)=\langle p, y+z\rangle +\langle x, q\rangle+ L(z,q)
\]
is convex and lower semi-continuous in the first variable $(z,q)$ and concave and upper semi-continuous in the second variable $p$, hence in view of the coercivity condition, Von-Neuman's min-max theorem applies and we get:
\begin{eqnarray*}
H(x,y)&=&\sup\limits_{p\in X^{*}}\inf\limits _{(z,q)\in X\times X^{*}}\{\langle p, y+z\rangle  +\langle x, q\rangle+ L(z,q)\}\\
&=&\inf\limits _{(z,q)\in X\times X^{*}}\sup\limits_{p\in X^{*}}\{\langle p, y+z\rangle  +\langle x, q\rangle+ L(z,q)\}\\
&=&\inf\limits _{q\in X^{*}}\{\langle x, q\rangle + L(-y,q)\}\\
&=&-H_L(-y,-x). 
\end{eqnarray*}
It follows that  $L$ is tempered  under the coercivity assumption. 

Suppose now $L(x,p_0) \leq C(1+\|x\|^\alpha)$, and  consider the $\lambda$-regularization of its conjugate $L^{*}$, that is 
$M_{\lambda}=L^{*}\star T^{*}_{\lambda}$ where $T^*_{\lambda}(x, p)=\frac{\|x\|^{2}}{2\lambda^{2}}+\frac{\lambda^{2}\|p\|^{2}}{2}$. Since obviously $L^{*}$ is ASD on $X^*$, we get from Lemma 2.5.3 that  $M_{\lambda}$ is ASD on $X^*$. Moreover, 
 \[
 M_{\lambda} (p,x) =\inf\{L^*(q,x)+ \frac{\|p-q\|^{2}}{2\lambda^{2}}+\frac{\lambda^{2}\|x\|^{2}}{2}; q\in X^*\}
\leq L(x, p_0)+\frac{1}{2\lambda^2}\|p-p_0\|^2 +\frac{\lambda^{2}\|x\|^{2}}{2}\leq C_1+C_2\|x\|^{\beta}+ C_3\|p\|^{2}
 \]
 which means that its dual $M_{\lambda}^{*}$ is an ASD Lagrangian on $X$ that is coercive in both variables. By the first part of the proof, $H_{M_{\lambda}^{*}}$ is a tempered ASD Hamiltonian on $X$.  
But in view of Lemma 2.5.1, we have  $M^{*}_{\lambda}=L\oplus T_{\lambda}$ and therefore $H_{M_{\lambda}^{*}}=H_{L}+H_{T_{\lambda}}$. Consequently $x\to H_{L}(x,y)$ is upper semi-continuous and $L$ itself is a tempered ASD Lagrangian. \\

 By exploiting  the duality between tempered ASD Lagrangians and ASD Hamiltonians, we get  the following   
 \begin{proposition} The  classe ${\cal L}_{_{TAD}}(X)$ possess the following permanence properties.
 \begin{enumerate}
\item  If $L$ and $M$ are in  ${\cal L}_{_{TAD}}(X)$ and $\lambda >0$,  then the Lagrangians $L\oplus M$,  and $\lambda {\bf \cdot} L$  also belong to  ${\cal L}_{_{TAD}}(X)$.  
\item  If $L$ is an ASD Lagrangian, then its $\lambda$-regularization $L_\lambda \in {\cal L}_{_{TAD}}(X)$. 
      \item If $L\in  {\cal L}_{_{TAD}}(X)$ and $\Lambda :X\to X^{*}$  is a skew-adjoint operator, then  $L_{ B}$ is also in  ${\cal L}_{_{TAD}}(X)$.
 
\end{enumerate}
\end{proposition}
\noindent{\bf Proof:} They all  follow from Proposition 2.1, Lemma 2.5 and Proposition 2.3. Note also that 
\[
H_{L_ B}(x,y)=H_L(x,y) -\langle  B x, y\rangle.
\]

Let now $B$ be a linear --not necessarily bounded- map from its domain $D(B)\subset X$ into $X^*$ such that  $D(B)$ dense in $X$, we consider the domain of its adjoint $B^*$ which is defined as:
\[
D(B^*)=\{x\in X; \sup \{\langle x,By\rangle; y\in D(B), \|y\|_X\leq 1\}<+\infty\}.
\]
\begin{definition}  Say that
\begin{enumerate}
\item  $B$ is antisymmetric if $D(B)\subset D(B^*)$ and if $B^*=-B$ on $D(B)$.
\item  $B$ is skew-adjoint if it is antisymmetric and if $D(B)= D(B^*)$. 
\end{enumerate}
\end{definition} 
We then have the  following easy lemma (See also \cite{GT2}).
\begin{lemma} Let $L:X\times X^*\to \R$ be an ASD Lagrangian on a reflexive Banach space $X$ and let $B$ be a linear  skew-adjoint map from its domain $D(B)\subset X$ into $X^*$ such that  the function $x\to L(x,0)$ is bounded on the unit ball  of $X$. 
The Lagrangian $L_B$  defined by
\begin{eqnarray*}
L_B(x,p)=\left\{ \begin{array}{l}
L(x,Bx+p) \quad 
\mbox{   if }x\in D(B)\\
+\infty 
\, \mbox{\phantom{XXXXXXX}if }x\notin D(B) 
\end{array}\right.
\end{eqnarray*}
is then itself anti-selfdual on $X$.  Moreover, if $L$ is tempered then so is $L_B$ whose Hamiltonian is given by  
\begin{eqnarray*}
H_B(x,y)=\left\{ \begin{array}{l}
H_L(x,y)-\langle Bx,y\rangle  \quad 
\mbox{   if }x\in D(B)\\
-\infty 
\, \mbox{\phantom{XXXXXXXXXX} if }x\notin D(B) 
\end{array}\right.
\end{eqnarray*}

\end{lemma}  
We shall also deal with situations where operators are skew-adjoint provided one takes into account certain boundary terms. We consider the following notion introduced in \cite{GT2}.

 \begin{definition}  \rm 
Let $B$ be a linear map from its domain $D(B)$ in a reflexive Banach space $X$ into $X^*$ and consider $(b_1, b_2)$ to be a pair of linear maps from its domain $D(b_1, b_2)$ in $X$ into the product of two Hilbert spaces $H_1 \times H_2$. Associate the set
\[
D^*(B,b_1, b_2)=\left\{y\in X;\, \sup\{\braket{y}{Bx}- \frac{1}{2}(\| b_1(x)\|_{H_1}^2+ \| b_2(x)\|_{H_2}^2); x\in S, \| x\|_X<1\} <\infty\right\}.
\]

$\bullet$\, Say that $B$ is {\it anti-symmetric modulo the boundary operators $(b_1,b_2)$} if the following properties are satisfied: 
 \begin{enumerate}
\item The set $S=D(B)\cap D(b_1, b_2)$ is dense in $X$.
\item  The space $X_0:=Ker (b_1,b_2)\cap D(B)$ is dense in $X$.
\item The image of $S$ by $(b_1,b_2)$ is dense in $H_1\times H_2$. 
\item For every $x, y \in S$, we have 
 $
\braket {y}{Bx}= -\braket{By}{x}+ \braket {b_2(x)}{b_2(y)}_{H_1} -\braket {b_1(x)}{b_1(y)}_{H_2}.
$
\end{enumerate}
$\bullet$ Say that $B$ is {\it skew-adjoint modulo the boundary operators $(b_1,b_2)$} if it is {\it anti-symmetric modulo the boundary operators $(b_1,b_2)$} and if in addition $D^*(B,b_1, b_2)=D(B)\cap D(b_1, b_2)$.
\end{definition}
It is clear that if $b_1$, $b_2$ are identically zero,  then our definition coincides with the notions  in Definition 3.5.

 For problems involving boundaries, we may start with an ASD Lagrangian $L$, but if linear operator $B$ is skew-adjoint modulo a term involving the boundary, the Lagrangian $L_B$ is not ASD  but we may recover anti-selfduality by adding a correcting term via a ``Boundary Lagrangian" $\ell$.
   \begin{definition} We say that $\ell :H_{1}\times H_{2} \to \R \cup \{+\infty\}$ is a {\it self-dual boundary Lagrangian} if
\begin{equation}
\ell^*(-h_{1},h_{2})=\ell (h_{1},h_{2}) \quad \hbox{\rm for all $(h_{1},h_{2}) \in H_{1}\times H_{2}$}.
\end{equation} 
\end{definition}
It is easy to see that such a boundary Lagrangian will always satisfy the inequality
\begin{equation}
\label{little.positivity}
\hbox{$\ell (r,s)\geq \frac{1}{2}(\|s\|^2-\|r\|^2)$ for all $(r,s)\in H_1\times H_2$.}
\end{equation}
The basic example of a self dual boundary Lagrangian is given by a function $\ell$ on $H_{1}\times H_{2}$, of the form $\ell (x,p) =\psi_{1}(x) +\psi_{2}(p)$, with  $\psi_{1}^{*}(x)=\psi_{1}(-x)$ and $ \psi_{2}^{*}(p)=\psi_{2}(p)$. Here the choices for $\psi_{1}$ and $\psi_{2}$ are rather limited and the typical sample is:
\[
\psi_{1} (x)= \frac{1}{2}\|x\|^{2} -2\langle a, x\rangle +\|a\|^{2}, \quad \hbox{\rm and \quad $\psi_{2}(p)=\frac{1}{2}\|p\|^{2}$.}
\]
where  $a$ is given in $H_{1}$. Boundary operators and Lagrangians allow us to build new ASD Lagrangians.
Here is the situation when the skew-adjoint operators are not necessarily bounded. Most of it  was established in \cite{GT2}, but we include here a proof for completness. 

 \begin{proposition} Let $\ell :H_{1}\times H_{2} \to \R \cup \{+\infty\}$ be a self-dual boundary Lagrangian on the product of two Hilbert spaces $H_{1}\times H_{2}$,  and let $L:X\times X^*\to \R$ be an ASD Lagrangian on a reflexive Banach space $X$ such that for ever $p\in X^*$, the function $x\to L(x,p)$ is bounded on the bounded sets of $X$. Let $B$ be a linear  map from its domain $D(B)\subset X$ into $X^*$, and let $(b_1,b_2):D(b_1,b_2)\subset X\to H_1\times H_2$ be linear boundary operators. Assume one of the following two conditions:
\begin{enumerate}
\item  $B$ is antisymmetric modulo the boundary operators $(b_1,b_2)$,  and $0\in {\rm Dom}_1(L) \subset D(B)\cap D(b_1, b_2)$.
\item  $B$ is skew-adjoint   modulo the boundary operators $(b_1,b_2)$ and $\ell (r,s) \leq C(1+\|r\|^2+\|s\|^2)$ for all $(r,s)\in H_1\times H_2$. 
\end{enumerate}
Then the Lagrangian defined by
 \begin{eqnarray*}
L_{B,\ell}(x,p)=\left\{ \begin{array}{l}
L(x,Bx+p)+ \ell (b_{1}(x), b_{2}(x))\quad 
\mbox{   if }x\in D(B)\cap D(b_1, b_2)\\
+\infty 
\, \mbox{\phantom{XXXXXXXXXXXXXXXXX}if }x\notin D(B)\cap D(b_1, b_2)
\end{array}\right.
\end{eqnarray*}
is anti-self dual on $X$. Its Hamiltonian is then given by
  \begin{eqnarray*}
H_{B,\ell}(x,y)=\left\{ \begin{array}{l}
H_L(x,y)-\langle Bx,y\rangle -  \ell (b_{1}(x), b_{2}(x))\quad 
\mbox{   if }x\in D(B)\cap D(b_1, b_2)\\
-\infty 
\, \mbox{\phantom{XXXXXXXXXXXXXXXXXXXX}if }x\notin D(B)\cap D(b_1, b_2)
\end{array}\right.
\end{eqnarray*}

\end{proposition}

\noindent {\bf Proof:}  Before we proceed we the proof, we note that while $H_{B,\ell}(-x,-y)\leq -H_{B,\ell}(y,x)$ for every $(x,y) \in X\times X$ and consequently $H_{B,\ell}(x,-x)\leq 0$, we almost never have equality unless it is zero on the boundary. In other words $L_{B,\ell}$ is never tempered even when $L$ is. The is due to the fact that with the above assumption on the density of their kernel, the operator $b_1, b_2$ can never be continuous. \\

 Assume now that $B$ is antisymmetric modulo the boundary operators $(b_1,b_2)$, and that for ever $p\in X^*$, the function $x\to L(x,p)$ is  continuous on $X$. We shall prove that $L_{B,\ell}^*(p,\tilde x)=L_{B,\ell}(-\tilde x,-\tilde p)$ if $\tilde x\in D(B)\cap D(b_1, b_2)$. Indeed, 
 fix  $\tilde x\in S:= D(B)\cap D(b_1, b_2)$, and write
\begin{eqnarray*}
L_{B,\ell}^*(\tilde p,\tilde x) &=&
\sup \left\{
   \braket{\tilde x}{p}+\braket{x}{\tilde p}
   -L(x,Bx+p)-\ell (b_{1}(x), b_{2}(x)); x\in S, p\in X^*\right\} 
\end{eqnarray*}
Substituting $q=Bx+p$, and since for $\tilde x\in S$,  we  have
$\braket{\tilde x}{Bx}=-\braket{x}{B\tilde x}-\braket{b_1(x)}{b_1(\tilde x)}
  +\braket{b_2(x)}{b_2(\tilde x)}$, 
we obtain
\begin{eqnarray*}
L_{B,\ell}^*(\tilde p,\tilde x) &=&
\sup\limits_{\stackrel{x\in S}{q\in X^*}} \left\{
  \braket{\tilde x}{q-Bx}+\braket{x}{\tilde p} 
   -L(x,q)-\ell (b_{1}(x), b_{2}(x))\right\}\\
 &=& \sup\limits_{\stackrel{x\in S}{q\in X^*}}
   \biggl\{ \braket{x}{B\tilde x}+\braket{b_1(x)}{b_1(\tilde x)}
   -\braket{b_2(x)}{b_2(\tilde x)} +\braket{\tilde x}{q}
    +\braket{x}{\tilde p} -L(x,q)-\ell (b_{1}(x), b_{2}(x))\biggl\}\\
 &=&
\sup\biggl\{\braket{x}{B\tilde x+\tilde p} +\braket{b_1(x+x_0)}{b_1(\tilde x)}
   -\braket{b_2(x+x_0)}{b_2(\tilde x)} +\braket{\tilde x}{q}-L(x,q) \\
 & &\quad \quad \quad \quad \quad -\ell (b_{1}(x+x_0), b_{2}(x+x_0));
    x\in S, q\in X^*,x_0\in\mbox{\rm Ker}(b_1,b_2)\cap D(B)\biggl\}
\end{eqnarray*}
Since $S$ is a linear space, we may set $w=x+x_0$ and write
 \begin{eqnarray*}
L_{B,\ell}^*(\tilde p,\tilde x) &=&
\sup\biggl\{\braket{w-x_0}{B\tilde x+\tilde p} 
   +\braket{b_1(w)}{b_1(\tilde x)} 
   -\braket{b_2(w)}{b_2(\tilde x)} 
   +\braket{\tilde x}{q} -L(w-x_0,q)\\
& &\quad \quad \quad -\ell (b_{1}(w), b_{2}(w));  w\in S, q\in X^*,x_0\in\mbox{Ker }(b_1,b_2)\cap D(B)\biggr\}
\end{eqnarray*}
Now, for each fixed $w\in S$ and $q\in X^*$, the supremum over $x_0\in$ ${\rm Ker} (b_1,b_2)\cap D(B)$ can be taken as a supremum over $x_0\in X$ since ${\rm Ker}(b_1,b_2)\cap D(B)$ is  dense in $X$ and all 
 terms involving $x_0$ are continuous in that variable. 
Furthermore, for each fixed $w\in S$ and $q\in X^*$, the supremum over $x_0\in X$ of the terms $w-x_0$ can be written as supremum over $v\in X$ where $v=w-x_0$. So setting $v=w-x_0$ we get
\begin{eqnarray*}
L_{B,\ell}^*(\tilde p,\tilde x) &=&
\sup\biggl\{ \braket{v}{B\tilde x+\tilde p} 
   +\braket{b_1(w)}{b_1(\tilde x)} 
   -\braket{b_2(w)}{b_2(\tilde x)} 
   +\braket{\tilde x}{q}-L(v,q)\\
& &\quad \quad \quad \quad \quad -\left.\ell (b_{1}(w), b_{2}(w)); 
    v\in X, q\in X^*, w\in S\right\}\\
&=&   \sup_{v\in X}\sup_{q\in X^*}\left\{\braket{v}{B\tilde x +\tilde p}
   +\braket{\tilde x}{q}-L(v,q) \right\}\\
 & &\quad\quad +\sup_{w\in S}\left\{\braket{b_1(w)}{b_1(\tilde x)} 
   +\braket{b_2(w)}{-b_2(\tilde x)} -\ell (b_{1}(w), b_{2}(w))\right\}
\end{eqnarray*}
Since the range of $(b_1,b_2):S\to H_1\times H_2$ is dense 
in the $H_1\times H_2$ topology, the boundary term can be written as
\[
  \sup_{a\in H_1}\sup_{b\in H_2}\left\{ 
   \braket{a}{b_1(\tilde x)}+\braket{b}{-b_2(\tilde x)}
   - \ell (a,b)\right\}
= \ell^* (b_1(\tilde x), -b_2(\tilde x))=\ell(-b_1(\tilde x), -b_2(\tilde x))..
  \]
while the main term is clearly equal to $ L^*(B\tilde x +\tilde p,\tilde x)=L(-\tilde x,-B\tilde x-\tilde p)$ in such a way that $L_{B,\ell}^*(p,\tilde x)=L_{B,\ell}(-\tilde x,-\tilde p)$ if $\tilde x\in D(B)\cap D(b_1, b_2)$. \\
  
  Now assume $\tilde x\notin S=D(B)\cap D(b_1, b_2)$, then $-\tilde x\notin S$. and we  distinguish the two cases:\\
  \noindent {\bf Case 1:} Under condition 1,   we have that  $-\tilde x \notin  {\rm Dom}_1(L)$, hence $-H_L(-\tilde x,0)=+\infty$. Now the boundedness condition on $L$ implies by Proposition 2.3 that it is a tempered ASD Lagrangian, which means since $0\in  {\rm Dom}_1(L)$ that  $H_L(0,\tilde x)= -H_L(-\tilde x,0).$ 
It follows
  \begin{eqnarray*}
  L_{B,\ell}^*(\tilde p,\tilde x) &=&
   \sup_{\stackrel{x\in  {\rm Dom}_1(L)}{p\in X^*}}\left\{
   \braket{\tilde x}{q-Bx}+\braket{x}{\tilde p} 
   -L(x,q)-\frac{{\| b_1(x)\|}_{H_1}^2}{2}
   -\frac{{\| b_2(x)\|}_{H_2}^2}{2} \right\} \\
& =& \sup_{p\in X^*}\left\{
   \braket{\tilde x}{p} 
   -L(0,p) \right\} \\
  &=&H_L(0,\tilde x)= -H_L(-\tilde x,0)=+\infty=L_{B,\ell}(-\tilde x, -\tilde p)
\end{eqnarray*}
  \noindent {\bf Case 2:} Under condition 2, write   
  \begin{eqnarray*}
L_{B,\ell}^*(\tilde p,\tilde x) &=&
   \sup_{\stackrel{x\in S}{q\in X^*}}\left\{
   \braket{\tilde x}{q-Bx}+\braket{x}{\tilde p} 
   -L(x,q)-\frac{{\| b_1(x)\|}_{H_1}^2}{2}
   -\frac{{\| b_2(x)\|}_{H_2}^2}{2} \right\} \\
& \geq&\sup_{\stackrel{x\in S}{{\| x\|}_X<1}}\left\{
   \braket{-\tilde x}{Bx}+\braket{x}{\tilde p} 
   -L(x,0)-\frac{{\| b_1(x)\|}_{H_2}^2}{2} 
   -\frac{{\| b_2(x)\|}_{H_2}^2}{2}.\right\} 
\end{eqnarray*}
Since by assumption $L(x,0)<K$ whenever ${\| x\|}_X<1$,  and $\ell (r,s) \leq C(1+\|r\|^2+\|s\|^2)$ for all $(r,s)\in H_1\times H_2$, we finally obtain that 
\begin{eqnarray*}
L_{B,\ell}^*(\tilde p,\tilde x)&\ge&\sup_{\stackrel{x\in S}{{\| x\|}_X<1}}\left\{
   \braket{-\tilde x}{Bx}+\braket{x}{\tilde p}-C-K
   -\frac{{\| b(x)\|}_{H_2}^2}{2} -\frac{{\| b_2(x)\|}_{H_2}^2}{2}\right\}\\
  &=& +\infty = L_{B,\ell}(-\tilde x,-\tilde p) 
  \end{eqnarray*}
since $\tilde x\notin S$ as soon as $-\tilde x\notin S$.
Therefore $ L_{B,\ell}^*(\tilde p,\tilde x)=L_{B,\ell}(-\tilde x,-\tilde p)$ for all $(\tilde x,\tilde p)\in X\times X^*$ and $L_{B,\ell}$ is an anti-selfdual Lagrangian. \\

\section{A nonlinear variational principle for ASD Lagrangians}

 \begin{definition} Say that a --non necessarily linear-- map $\Lambda:D(\Lambda)\subset X\to X^{*}$ is a {\it regular conservative map} if it is weak-to-weak continuous and if   $\langle \Lambda x, x\rangle =0$ for all $x$ in its domain  $D(\Lambda)$.
 \end{definition}
 It is clear that  skew-symmetric bounded linear operators are regular conservative maps. However, there are also plenty of nonlinear ones  many of them appearing in the basic equations of  hydrodynamics and magnetohydrodynamics (see below and \cite{Te}).    

Now recall  that ASD Lagrangians readily satisfy
$L(x, p)\geq -\langle x, p\rangle$ for every $(x, p) \in X\times X^{*}$,  
which means that they are non-negative on the graphs of conservative maps, that is: 
\begin{equation}
\label{always.positive}
L(x,\Lambda x)\geq 0 \quad \hbox{\rm for all $x\in D(\Lambda)$}.
\end{equation}
What is remarkable is that, just like in the case of linear skew-adjoint operators \cite{G2}, the infimum will often be zero, a fact that will allow us to derive variationally several nonlinear  PDEs  without using Euler-Lagrange theory. Here is our basic result.

\begin{theorem} Let $\ell :H_{1}\times H_{2} \to \R \cup \{+\infty\}$ be a self-dual boundary Lagrangian on the product of two Hilbert spaces $H_{1}\times H_{2}$, and let $L:X\times X^{*}\to {\bf R}\cup\{+\infty\}$ be an  anti-selfdual Lagrangian on  a reflexive Banach space $X$ such that   its Hamiltonian $H_{L}$ satisfies  $\lim\limits_{\|x\|\to +\infty}H_{L}(0,x)=+\infty$.\\
 Let $B: D(B)\subset X \to X^*$ and $(b_{1}, b_{2}): D(b_1, b_2)\subset X\to  H_{1}\times H_{2}$ be linear operators such that 
 \begin{equation}
\hbox{$ \langle   x, B  x\rangle= \frac{1}{2}(\|b_{2}x\|^{2} -\|b_{1}x\|^{2}) $ for all $x\in  D(A)\cap D(b_1, b_2)$,}
 \end{equation}
   Suppose $\Lambda: D(\Lambda)\subset X\to X^{*}$  is a regular conservative operator such that the Lagrangian $L_{B,\ell}$ is anti-selfdual on the graph of $\Lambda$ and ${\rm Dom}_1(L)\cap D(B)\cap D(b_1, b_2) \subset D(\Lambda)$. \\
   Then, there exists $\bar x \in {\rm Dom}_1(L)\cap D(B)\cap D(b_1, b_2)$ such that:
  \begin{eqnarray}
\label{eqn:zero}
 L(\bar x, B \bar x +\Lambda \bar x) +  \ell (b_{1}\bar x, b_{2}\bar x)&=&\inf_{x\in X} \left\{L(x, Bx+\Lambda x) +  \ell (b_{1}x, b_{2}x)\right\}= 0\\
  \hfill (-\Lambda \bar x-B \bar x, -\bar x) &\in & \partial L (\bar x,B\bar x+\Lambda \bar x)\\
  \ell (b_1(\bar x), b_2(\bar x))&=& \frac{1}{2}(\|b_{2}(\bar x)\|^{2} -\|b_{1}(\bar x)\|^{2}). 
 \end{eqnarray}
  \end{theorem}
We shall deduce Theorem 3.2 from the following Ky-Fan type min-max theorem due to Brezis-Nirenberg-Stampachia (see \cite{BNS}).
 \begin{lemma} Let $D$ be a  convex subset of a reflexive Banach space $X$ and let $M(x,y)$ be a real valued function  on $D\times D \subset X\times X$ that satisfies the following conditions:
\begin{description}
\item (1) $M(x,x) \leq 0$ for every $x\in D$.
\item (2) For each $x\in D$,  the function  $y \to M(x,y)$ is concave.
\item (3) For each $y\in D$, the function $x\to M(x,y)$ is weakly lower semi-continuous on $X$.
\item (4) There exists $K>0$ and $y_{0}\in X$ such that $\|y_{0}\|\leq K$ and $\inf\limits_{\|x\|>K}M(x,y_{0})>0$. 
\end{description} 
Then there exists $x_{0}\in D$ such that $\sup\limits_{y\in D}M(x_{0},y)\leq 0$.
\end{lemma}
\noindent{\bf Proof of Theorem 3.2:} Since the Lagrangian  $L_{B,\ell} $ defined above 
  is anti-self dual on the graph of $\Lambda$, we can write for each $x\in D:={\rm Dom}(L)\cap D(B)\cap D(b_1, b_2)\subset D(\Lambda)$,  
  \begin{eqnarray*}
  I(x)&=&L_{B,\ell}  (x, \Lambda x)=L_{B,\ell} ^{*}(-\Lambda x, -x)\\
&=&\sup\{\langle y,-\Lambda x\rangle +\langle p, -x\rangle -L_{B,\ell} (y,p); y\in X, p\in X^{*}\}\\
&=&\sup\{\langle y,-\Lambda x\rangle +\langle p, -x\rangle - L(y, B y+p) -\ell (b_{1}(y), b_{2}(y)); y\in D, p\in X^{*}\}\\
&=&\sup\{\langle y,-\Lambda x\rangle +\langle q-By, -x\rangle - L(y, q) -\ell (b_{1}(y), b_{2}(y)); y\in D, q\in X^{*}\}\\
&=&\sup\{\langle y,-\Lambda x\rangle +\langle x,By\rangle-\ell (b_{1}(y), b_{2}(y))+\sup\{\langle q, -x\rangle - L(y, q); , q\in X^{*}\}; y\in D\}\\
&=&\sup\{\langle y,-\Lambda x\rangle +\langle x,By\rangle-\ell (b_{1}(y), b_{2}(y))+H_{L}(y, -x); \, y\in D\}\\
&=& \sup\limits_{y\in D}M(x,y) 
\end{eqnarray*}
where $
M(x,y)= \langle y,-\Lambda x\rangle +\langle x,By\rangle-\ell (b_{1}(y), b_{2}(y))+H_{L}(y, -x)$, 
 and 
where  $H_L$ is the Hamiltonian associated to $L$. \\
We now claim that $M$ satisfies all the properties of the Ky-Fan min-max lemma above. Indeed, 

(1)  For each $x\in D(B)\cap D(b_1,b_2)$, we have  $y \to M(x,y)$ is concave since the first part $y\to \langle y,-\Lambda x\rangle+\langle x,By\rangle$ is clearly linear, while $y\to -\ell (b_{1}(y), b_{2}(y))$ and $y\to H_L(y,x)$ are concave.

(2) For each $y\in D(A)\cap D(b_1, b_2)$, the function $x\to M(x,y)$ is weakly lower semi-continuous on $D(A)\cap D(b_1, b_2)$ since 
$x\to \langle y,-\Lambda x\rangle+\langle x,By\rangle$ is weakly continuous  while $x\to H_L(y, -x)$ is clearly the supremum of continuous affine functions. 

(3) To show that $M(x,x) \leq  0$ for each $x\in D(A)\cap D(b_1, b_2)$, use the fact that the operator $\Lambda$ is conservative (i.e.,  $\langle \Lambda x, x \rangle =0$ on $D(\Lambda)$) and that $H_L$ is an ASD Hamiltonian, hence $H_L(x,-x) \leq 0$ to write
\[
M(x,x) \leq \langle x,Bx\rangle-\ell (b_{1}(x), b_{2}(x))=\frac{1}{2}(\|b_2(x)\|^2-\|b_1(x)\|^2)-\ell (b_{1}(x), b_{2}(x))\leq 0.
\]
(4) The set $X_{0}=\{x\in X; M(x, 0)\leq 0\}$ is bounded in $X$ since $M(x,0)=H_L(0,-x)-\ell (0,0)$ and the latter goes to infinity with $\|x\|$. 

It follows from Lemma 3.3 that there exists $\bar x\in D$ such that $\sup\limits_{y\in D}M(\bar x,y)\leq 0$. In other words
\[
I(\bar x)= \sup\limits_{y\in D}M(\bar x,y) \leq 0.
\]
On the other hand, for any $x\in X$, we have
\begin{eqnarray*}
I(x)&=& L(\bar x, B\bar x+ \Lambda \bar x) +  \ell (b_{1}\bar x, b_{2}\bar x)\\
&\geq& - \langle  x, B x+ \Lambda  x\rangle +  \ell (b_{1} x, b_{2} x)\\
&=&- \langle  x, B x \rangle +  \ell (b_{1} x, b_{2} x)\\
&=&- \frac{1}{2}(\|b_2(x)\|^2-\|b_1(x)\|^2) +  \ell (b_{1} x, b_{2} x)\geq 0.
\end{eqnarray*}
 It follows that $I(\bar x)=0=\inf_{x\in X}I(x)$, which means 
 \begin{equation}
\label{eqn:zero.2}
 L(\bar x, B\bar x+ \Lambda \bar x) +  \ell (b_{1}\bar x, b_{2}\bar x)=\inf_{x\in X} \left\{L(x, Bx+\Lambda x) +  \ell (b_{1}x, b_{2}x)\right\}= 0.
 \end{equation}
 To establish (\ref{eqn:zero}), write 
 \begin{eqnarray*}
 0= L(\bar x, B \bar x+ \Lambda \bar x) +  \ell (b_{1}\bar x, b_{2} \bar x)&=& L( \bar x, B \bar x+\Lambda \bar x)+ \langle  \bar x, B\bar x+\Lambda \bar x\rangle- \langle \bar x, B\bar x+\Lambda  \bar x\rangle + \ell (b_{1} \bar x, b_{2} \bar x)\\
  &=&L(\bar x, B \bar x+\Lambda \bar x)+ \langle \bar x, \Lambda \bar x\rangle- \frac{1}{2}(\|b_{2}\bar x\|^{2} -\|b_{1}\bar x\|^{2})  + \ell (b_{1} \bar x, b_{2} \bar x). 
\end{eqnarray*}  
  Since $L(x,p)+\langle x,p\rangle \geq 0$ and $\ell (r,s)\geq \frac{1}{2}(\|s\|^2-\|r\|^2)$, we get
   \begin{equation}
 \left\{ \begin{array}{lcl}
\label{eqn:application}
  L(\bar x, B\bar x+\Lambda \bar x) + \langle \bar x, B \bar x\rangle
 &=&0.\\
\hfill  \ell (b_{1} \bar x, b_{2} \bar x)&=&\frac{1}{2}(\|b_{2}\bar x\|^{2} -\|b_{1}\bar x\|^{2}).
 \end{array}\right.
 \end{equation}
 To obtain the second claim, we use that $L$ is anti-selfdual to write
\begin{eqnarray*}
\langle (\bar x, \Lambda \bar x+B \bar x), (-\Lambda \bar x-B \bar x, -\bar x)\rangle&=&0\\
&=& 2L( \bar x, \Lambda \bar x+B \bar x)\\
&=&L( \bar x, \Lambda \bar x+B \bar x)+L^{*}(-\Lambda \bar x-B \bar x, -\bar x).
\end{eqnarray*}
The last part of claim (\ref{eqn:zero}) now follows from the limiting case of the Legendre-Fenchel duality.

\begin{remark} \rm The above holds under an asssumption like 
$
L(0,p)\leq C\|p\|^{\alpha}
$
where $\alpha>1$, since then $H_L(0,-x)$ is coercive being its Legendre dual in the $p$-variable and therefore $M(x,0)$ is coercive.
 \end{remark} 
 The following corollary is immediate as it corresponds to the case where $b_1, b_2$ and $\ell$ are zero, while $B:D(B)\subset X\to X^{*}$ is  any linear  antisymmetric  operator with a large enough domain.  We shall see in the next section that it is already sufficient to cover several nonlinear PDEs including Navier-Stokes equations and others.  

\begin{corollary}  Let $L$ be an anti-self dual Lagrangian on a reflexive  space $X$ and let $H_{L}$ be  the corresponding Hamiltonian such that $\lim\limits_{\|x\|\to +\infty}H_{L}(0,x)=+\infty$. Let $B$ be an antisymmetric linear operator such that  ${\rm Dom}_1(L) \subset D(B)$. 
Then for any regular conservative operator $\Lambda:D(\Lambda)\subset X\to X^{*}$ such that $ {\rm Dom}_1( L)\subset D(\Lambda)$,  there exists $\bar x \in {\rm Dom}_1(L)$ such that:
  \begin{eqnarray}
 L(\bar x, B \bar x +\Lambda \bar x) &=&\inf_{x\in X} L(x, Bx+\Lambda x) = 0\label{eqn:zero.1}\\
  \hfill (-\Lambda \bar x-B \bar x, -\bar x) &\in & \partial L (\bar x,B\bar x+\Lambda \bar x).\label{eqn:zero.2}
  \end{eqnarray}
\end{corollary}
\noindent{\bf Proof:} It is sufficient to apply Theorem 3.2 to $b_1, b_2$ and $\ell$ being identically zero, while $\tilde \Lambda=\Lambda +B$  satisfies $ {\rm Dom}_1 (L)\subset D(\tilde \Lambda)$. \\

If the domain of the anti-symmetric operator $B$ is not large enough, we can use Lemma 2.7 to obtain
  
\begin{corollary} Let $L$ be an anti-self dual Lagrangian on a reflexive Banach space $X$ such that  $\lim\limits_{\|x\|\to +\infty}H_{L}(0,x)=+\infty$, where $H_L$ is the corresponding Hamiltonian.  Let $B:D(B)\subset X\to X^{*}$ be a  linear  skew-adjoint operator such that $x\to L(x,0)$ is bounded on the unit ball of $X$.
 Then, for  any regular conservative operator $\Lambda :X\to X^{*}$ such that $D(B)\cap  {\rm Dom}_1 L\subset D(\Lambda)$, there exists $\bar x\in D(B)\cap  {\rm Dom}_1 L$ satisfying (\ref{eqn:zero.1}) and (\ref{eqn:zero.2}). 
\end{corollary}
{\bf Proof:} By the above Lemma, $L_B$ is an anti-selfdual Lagrangian, in particular it is so on the graph of $\Lambda$.  The rest follows from Theorem 3.2.\\

In order to deal with situations where operators are skew-adjoint provided one takes into account certain boundary terms, we have the following
   
  \begin{corollary} Let $\ell :H_{1}\times H_{2} \to \R \cup \{+\infty\}$ be a self-dual boundary Lagrangian on the product of two Hilbert spaces $H_{1}\times H_{2}$, and let $L:X\times X^{*}\to {\bf R}\cup\{+\infty\}$ be an  anti-selfdual Lagrangian on  a reflexive Banach space $X$ such that 
  \[
  \hbox{$\lim\limits_{\|x\|\to +\infty}H_{L}(0,x)=+\infty$ and for ever $p\in X^*$,  $x\to L(x,p)$ is bounded on the bounded sets of $X$.}
  \]
  \\
   Let $(B, b_1,b_2): X\to H_1\times H_2$ be linear  operators such that one of the following two conditions hold:
\begin{enumerate}
\item  $B$ is antisymmetric modulo boundary operators $(b_1,b_2)$,  and $0\in {\rm Dom}_1(L) \subset D(B)\cap D(b_1, b_2)$.
\item  $B$ is skew-adjoint   modulo the boundary operators $(b_1,b_2)$ and $\ell (r,s) \leq C(1+\|r\|^2+\|s\|^2)$ for all $(r,s)\in H_1\times H_2$. 
\end{enumerate}
   Then for any regular conservative operator $\Lambda: D(\Lambda)\subset X\to X^{*}$ such that  ${\rm Dom}(L)\cap D(B)\cap D(b_1, b_2) \subset D(\Lambda)$, there exists $\bar x \in {\rm Dom}_1(L)\cap D(B)\cap D(b_1, b_2)$ such that:
  \begin{eqnarray}
\label{eqn:zero.12}
 L(\bar x, B \bar x +\Lambda \bar x) +  \ell (b_{1}\bar x, b_{2}\bar x)&=&\inf_{x\in X} \left\{L(x, Bx+\Lambda x) +  \ell (b_{1}x, b_{2}x)\right\}= 0\\
  \hfill (-\Lambda \bar x-B \bar x, -\bar x) &\in & \partial L (\bar x,B\bar x+\Lambda \bar x)\\
  \ell (b_1(\bar x), b_2(\bar x))&=& \frac{1}{2}(\|b_{2}(\bar x)\|^{2} -\|b_{1}(\bar x)\|^{2}). 
 \end{eqnarray}
  \end{corollary}
\noindent{\bf Proof:} This follows from Theorem 3.2 and Proposition 2.5.

 \section{A variational nonlinear Lax-Milgram theorem and applications}
  
  We now apply the above results to the most basic ASD Lagrangians of the form $L(x,p)=\phi(x)+\phi^*(Bx-p)$ where $\phi$ is a convex function and $B$ is a linear anti-symmetric but not necessarily bounded operator. The  applications differ as they  will depend on the ``position" of the domain of  $B$. We start with the case where the linear operator component has a ``large domain".
 
  \begin{theorem}  Let $\phi$ be a proper convex lower semi-continuous function on a reflexive Banach space $X$ such that $\lim\limits_{\|x\|\to +\infty}\frac{\phi(x)}{\|x\|}=+\infty$ and let $B:D(B)\subset X\to X^{*}$ be an anti-symmetric  linear operator such that ${\rm Dom}(\phi)\subset D(B)$. Then, for any $f\in X^{*}$ and any regular conservative operator $\Lambda: D(\Lambda)\subset X\to X^{*}$ such that ${\rm Dom}(\phi)\subset D(\Lambda)$, there exists a solution $\bar x \in {\rm Dom}(\phi)$ to the equation 
    \begin{equation}
 \label{Lax.Mil.0}
 0\in  f+\Lambda x+Bx+\partial \varphi (x).
  \end{equation}
  It is obtained as a minimizer of the problem:
\begin{equation}
 \label{min.99}
\inf_{x\in X}\left\{\phi (x) +\langle f, x\rangle +\phi^{*}(-\Lambda x- Bx-f) \right\}=0  
\end{equation}
   \end{theorem}
\noindent {\bf Proof:} It is an immediate consequence of Corollary 3.5 applied to the Lagrangian $L(x,p)=\psi(x)+\psi^*(-p)$ where $\psi (x)= \phi (x) +\langle f, x\rangle$.  Note that its Hamiltonian is now 
$
H(x,y)=- \langle Bx, y\rangle +\phi (-y)-\phi (x) -\langle f, x+y\rangle
$
meaning that the coercivity hypothesis  implies that $H(0,y)\to +\infty$ with $\|y\|$. Corollary 3.5 then applies with the Lagrangian $L$ and the conservative operator $\Lambda$ to obtain that the minimum in (\ref{min.99}) is attained at some $\bar x \in X$. We then get
\[
\phi (\bar x) +\phi^{*}(-B \bar x -\Lambda \bar x-f)=\langle -B \bar x-\Lambda \bar x-f ,\bar x\rangle
\]
 which yields, in view of Legendre-Fenchel duality that 
 \[
- B \bar x -\Lambda \bar x -f\in   \partial \varphi (\bar x). 
\] 
 An immediate application is the case where the linear operator component is bounded  which already covers many interesting applications.
  
  \begin{corollary}  Let $\phi$ be a  function on a reflexive Banach space $X$ and let $B:X\to X^{*}$ be a  bounded linear operator such that the function $\psi (x):=\phi(x)+\frac{1}{2} \langle Bx, x\rangle$ is  proper convex and lower semi-continuous on $X$. Assume
  \begin{equation}
  \label{coercive.function}
    \lim\limits_{\|x\|\to \infty}\|x\|^{-1}(\varphi (x) + \frac{1}{2}\langle Bx, x\rangle) =+\infty. 
 \end{equation}
  Then, for any regular conservative operator $\Lambda:X\to X^{*}$ and any $f\in X^{*}$, there exists a solution $\bar x \in X$ to the equation 
 \begin{equation}
 \label{Lax.Mil.1}
 0\in  f+\Lambda x+Bx+\partial \varphi (x).
  \end{equation}
  It is obtained as a minimizer of the problem:
\begin{equation}
 \label{min.100}
\inf_{x\in X}\left\{\psi (x) +\langle f, x\rangle  +\psi^{*}(-\Lambda x- B^{a}x-f) \right\}=0  
\end{equation}
where  $B^{a}$ is the anti-symmetric part of $B$.
   \end{corollary}
\noindent{\bf Proof:} Apply the above theorem to  $\psi (x) +\langle f, x\rangle $ and $B^a=\frac{1}{2}(B-B^*)$. The above theorem applies and we get   $\bar x \in X$ such that:
  \[
- B^{a}\bar x -\Lambda \bar x +f\in  \partial \psi (\bar x)= B^{s}\bar x  +\partial \varphi (\bar x) 
\] 
hence $\bar x$ satisfies (\ref{Lax.Mil.1}).\\

 \subsection*{Example 1: A variational resolution for the Stationary Navier-Stokes equation}
 Consider the incompressible stationary Navier-Stokes equation on a domain $\Omega$ of $\R^{3}$
 \begin{equation}
\label{NSE1}
 \left\{ \begin{array}{lcl}
    \hfill
 (u\cdot \nabla)u +f &=&\nu \Delta u - \nabla  p \quad \hbox{\rm on $ \Omega$}\\
\hfill {\rm div} u&=&0 \quad \hbox{\rm on  $\Omega$}\\
\hfill u &=&0 \quad \hbox{\rm on $\partial \Omega$}\\
\end{array}\right.
\end{equation}
where $\nu >0$ and $f\in L^{p}(\Omega;\R^{3})$. Let 
\begin{equation}
\label{phi}
\Phi (u)=\frac{\nu}{2} \int_{\Omega}\Sigma_{j,k=1}^{3}(\frac {\partial u_{j} } {\partial x_{k}})^{2}\, dx
\end{equation}
be the convex and coercive function on the Sobolev subspace $V=\{u\in H^{1}_0(\Omega; {\bf R}^{3}); {\rm div} v=0\}$. Its Legendre transform $\Phi^{*}$ on $V$  can be characterized as $  \Phi^{*} (v)= \langle Sv,v\rangle$ where $S:V^*\to V$ is the bounded linear operator that associates to $v\in V^*$ the solution $\hat v=Sv$ of the Stokes' problem
\begin{equation}
\label{Stokes}
 \left\{ \begin{array}{lcl}
    \hfill
  \nu \Delta \hat v + \nabla p&=&-\hat v \quad \hbox{\rm on $ \Omega$}\\
\hfill {\rm div} \hat v&=&0 \quad \hbox{\rm on  $\Omega$}\\
\hfill \hat v &=&0 \quad \hbox{\rm on $\partial \Omega$.}\\
\end{array}\right.
\end{equation}
It is easy to see that (\ref{NSE1}) can be reformulated as
 \begin{equation}
\label{NSE2}
 \left\{ \begin{array}{lcl}
    \hfill
 (u\cdot \nabla)u +f &\in& -\partial \Phi (u)=\nu \Delta u - \nabla  p  \\
\hfill  u&\in& V. \\
\end{array}\right. 
\end{equation}
Consider now the nonlinear operator $\Lambda: V \to V^{*}$ defined as
\[
\langle \Lambda u, v\rangle =\int_{\Omega}\Sigma_{j,k=1}^{3}u_{k}\frac {\partial u_{j} } {\partial x_{k}}v_{j}\, dx=\langle (u\cdot \nabla)u,v\rangle.
\]
We can deduce the following
\begin{theorem} Assume $\Omega$ is bounded domain in $\R^{3}$ and consider $f\in L^{p}(\Omega;\R^{3})$ for $p>\frac{6}{5}$. Then the infimum of the functional 
 \[
{I}(u)=\Phi (u)+\Phi^{*}(-(u\cdot \nabla)u+f )-\int_{\Omega}\Sigma_{j=1}^{3}f_{j}u_{j}
\]
 on $V$ is equal to zero, and is attained at a solution of the Navier-Stokes equation (\ref{NSE1}).
  \end{theorem}
\noindent {\bf Proof:} To apply Theorem  4.1, it remains to show that $\Lambda$ is a regular conservative operator. First note that $\langle \Lambda u, u \rangle =0$ on $V$. For the weak-to weak continuity,  assume that  $u^{n}\to u$ weakly in $H^{1}(\Omega)$. We need to show that for a  fixed $v\in V$, we have that 
\[
\langle \Lambda u^{n}, v\rangle=\int_{\Omega}\Sigma_{j,k=1}^{3}u^{n}_{k}\frac {\partial u_{j}^n } {\partial x_{k}}v_{j}\, dx=-\int_{\Omega}\Sigma_{j,k=1}^{3}u^{n}_{k}\frac {\partial v_{j}} {\partial x_{k}}u^n_{j}\, dx
\]
 converges  to  $\langle \Lambda u,v\rangle=\int_{\Omega}\Sigma_{j,k=1}^{3}u_k\frac {\partial v_{j} } {\partial x_{k}}u_{j}\, dx$. But the Sobolev embedding in dimension 3 implies that $(u^{n})$ converges strongly in $L^{p}(\Omega;\R^{3})$ for $1\leq p <6$. On the other hand, $\frac {\partial u_{j} } {\partial x_{k}}$ is in $L^{2}(\Omega)$ and the result follows from  an application of H\"older's inequality.
  
\subsection*{Example 2: Variational resolution for a fluid driven by its boundary}
The full strength of Corollary 4.2 comes out when one deals with the Navier-Stokes equation with  a boundary moving with a prescribed velocity:
\begin{equation}
\label{NSE3}
 \left\{ \begin{array}{lcl}
    \hfill
 (u\cdot \nabla)u +f &=&\nu \Delta u - \nabla  p \quad \hbox{\rm on $ \Omega$}\\
\hfill {\rm div} u&=&0 \quad \hbox{\rm on  $\Omega$}\\
\hfill u&=&u^0 \quad \hbox{\rm on $\partial \Omega$}\\
\end{array}\right.
\end{equation}
where $\int_{\partial \Omega} u^{0} {\bf \cdot n}\, d\sigma =0$, $\nu>0$ and $f\in L^{p}(\Omega;\R^{3})$.  Assuming  that $u^{0} \in H^{3/2}(\partial \Omega)$ and that $\partial \Omega$ is connected,  a classical result of Hopf then yields for each $\epsilon >0$, the existence of $v^0\in H^{2}(\Omega)$ such that
\begin{equation}
\label{Hopf}
v^0=u^{0}\,\, \hbox{\rm on $\partial \Omega$,\quad  ${\rm div}\, v^0=0$\quad  and \quad  $\int_{\Omega}\Sigma_{j,k=1}^{3}u_{k}\frac {\partial v^0_{j} } {\partial x_{k}}u_{j}\, dx\leq \epsilon \|u\|^{2}_{V}$ for all $u\in V$.}
\end{equation}
Setting $v=u+v^0$, then solving  (\ref{NSE3}) reduces to finding a solution for
\begin{equation}
\label{NSE4}
 \left\{ \begin{array}{lcl}
    \hfill
 (u\cdot \nabla)u +(v^0\cdot \nabla)u +(u\cdot \nabla)v^0 +f -\nu \Delta v^0 + (v^0\cdot \nabla)v^0 &=&\nu \Delta u - \nabla  p \quad \hbox{\rm on $ \Omega$}\\
\hfill {\rm div} u&=&0 \quad \hbox{\rm on  $\Omega$}\\
\hfill u&=&0 \quad \hbox{\rm on $\partial \Omega.$}\\
\end{array}\right.\nonumber 
\end{equation} 
This can be reformulated as the following equation in the space $V$
\begin{equation}
\label{NSE5}
  (u\cdot \nabla)u +(v^0\cdot \nabla)u +(u\cdot \nabla)v^0 +g   \in -\partial \Phi (u) 
\end{equation}
where $\Phi$ is again the convex functional $\Phi (u)=\frac{\nu}{2} \int_{\Omega}\Sigma_{j,k=1}^{3}(\frac {\partial u_{j} } {\partial x_{k}})^{2}\, dx$ as above and where 
\[
g:= f -\nu \Delta v^0 + (v^0\cdot \nabla)v^0\in V^{*}.
\]
In other words, this is an equation of the form
\begin{equation}
 \label{Lax.Mil.again}
  \Lambda u+Bu+g\in -\partial \Phi (u) 
  \end{equation}
  with $\Lambda u= (u\cdot \nabla)u$ is a regular conservative operator, and $Bu=(v^0\cdot \nabla)u +(u\cdot \nabla)v^0$ is a bounded linear operator. Note that the component $B^1u:=(v^0\cdot \nabla)u$ is skew-symmetric which means that  Hopf's result yields the required coercivity condition:
  \[
  \Psi (u):=\Phi (u)+\frac{1}{2}\langle Bu,u\rangle \geq \frac{1}{2}(\nu-\epsilon)\|u\|^{2} \quad {\rm for\,  all} \,\, u\in V.
  \]
  In other words,  $\Psi$  is convex and coercive and therefore we can apply Theorem 4.1 to deduce 
  \begin{theorem} Under the above hypothesis, and letting $A^a$ be the antisymmetric part of the operator $Au=(u\cdot \nabla)v^0$, the following functional 
  \[
   {I}(u)=\Psi (u)+\Psi^{*}(-(u\cdot \nabla)u-(v^0\cdot \nabla)u-A^au+g )-\int_{\Omega}\Sigma_{j=1}^{3}g_{j}u_{j}
\]
  has zero for infimum on the Banach space $V$, which is  attained at a solution $\bar u$ for (\ref{NSE4}).  
  \end{theorem}
  
  The next application is a nonlinear Lax-Milgram type result with boundary constraints.

\begin{theorem}
Let $\phi :X\to \R\cup \{+\infty\}$ be a convex and lower semi-continuous on a reflexive Banach space such that for some constant $C>0$ and $p_1,p_2>1$, we have
\begin{equation}
\hbox{$\frac{1}{C}\left(\| x\|_X^{p_1}-1\right)\leq\phi (x) \leq C\left( \| x\|_X^{p_2}+1\right)$    for every $x\in X$,}
\end{equation}
 Let  $B: D(B)\subset X \to X^*$ be a skew-adjoint operator modulo the boundary  $(b_1, b_2):D(b_1,b_2)\to H_1\times H_2$ where $H_1$, $H_2$ are two Hilbert spaces.
 Then for any regular conservative operator $\Lambda: X\to X^*$ such that $D(B)\cap D(b_1, b_2)\subset D(\Lambda)$ and any $a \in X$, there exists a solution $\bar x \in X$ to the equation 
 \begin{eqnarray}
 \label{Lax.Mil.2}
 \Lambda x+Bx+f&\in&- \partial \varphi (x)\\
 b_1(\bar x)&=&a.\nonumber
  \end{eqnarray}
  It is obtained as a minimizer of the functional defined as:
\begin{equation}
 \label{min.100}
 I(x)=\phi (x) + \langle f, x\rangle +\phi^{*}(-\Lambda x- Bx-f)  +\frac{1}{2}(\|b_{1}(x)\|^{2} +\|b_{2}(x)\|^{2} )-2\langle a, b_{1}(x)\rangle +\|a\|^{2} 
\end{equation}
when  $x\in D(A)\cap D(b_1, b_2)$ and $+\infty$ elsewhere. Moreover, $I(\bar x)=\inf_{x\in X}I(x)=0$. 
 \end{theorem}
 
\noindent {\bf Proof:} Let $\psi (x)=\phi (x) + \langle f, x\rangle$ and apply Corollary 3.7 to the ASD Lagrangian $L(x,p):=\psi (x)+\psi ^* (-p)$, to the boundary Lagrangian $\ell (r,s)= \frac{1}{2}(\|r\|^{2} +\|s\|^{2} )-2\langle a, r\rangle +\|a\|^{2} $, and to the skew-adjoint  triplet $(B, b_1,b_2)$. Note also that  $I$ can be rewritten as:
 \[
 I(x)= \phi (x) +\phi^{*}(-\Lambda x- Bx-f) -\langle x, -\Lambda x -Bx-f \rangle + \frac{1}{2}(\|b_{1}(x)-a\|^{2}.
    \]

  \subsection*{Example 3: Variational resolution for a fluid driven by a transport operator}
  
Let  $\vec a \in C^{\infty}(\bar\Omega)$ be a smooth vector field on a neighborhood of a $C^\infty$ bounded open set  $\Omega \subset \R^3$, let  $a_0 \in L^{\infty}(\Omega)$, and consider  the space $X = \{u\in H^1_0(\Omega; \R^3); {\rm div}(u)=0\}$ and the transport  operator $B: u\mapsto  a\cdot \nabla u+\frac{1}{2} {\rm div } u$ from $D(B)=\{u\in X; a\cdot \nabla u+\frac{1}{2} {\rm div } u \in X^*\}$ into $X^*$.  It is easy to show using Green's formula that the operator $b$ is skew-adjoint modulo the boundary $(b_1,b_2)$ on the space $X$ (See \cite{GT2}).  
Consider now the following equation on the domain $\Omega \subset \R^{3}$
 \begin{equation}
\label{NSE5}
 \left\{ \begin{array}{lcl}
    \hfill
 (u\cdot \nabla)u + \vec a\cdot  \nabla  u + a_0 u+|u|^{m-2}u+f &=&\nu \Delta u - \nabla  p \quad \hbox{\rm on $ \Omega$}\\
\hfill {\rm div} u&=&0 \quad \quad \quad \quad \quad \hbox{\rm on  $\Omega$}\\
\hfill u &=&0  \quad \quad \quad \quad \quad \hbox{\rm on $\partial \Omega$}\\
\end{array}\right.
\end{equation}
where $\nu >0$, $6\geq m\geq 1$  and $f\in L^{q}(\Omega;\R^{3})$ for $q\geq \frac{6}{5}$. Suppose 
\begin{equation}
\label{pos}
 \frac{1}{2}{\rm div}  (a) - a_0 \geq 0 \quad {\rm on} \quad \Omega, 
\end{equation} 
and consider the functional 
\begin{equation}
\label{phi}
\Psi (u)=\frac{\nu}{2} \int_{\Omega}\Sigma_{j,k=1}^{3}(\frac {\partial u_{j} } {\partial x_{k}})^{2}\, dx+\frac{1}{4}\int_\Omega( {\rm div} \vec a -2 a_0)\vert u\vert^2dx +\frac{1}{m}\int_\Omega \vert u\vert^mdx +\int_\Omega ufdx  
\end{equation}
which is  convex and coercive function on the space $X$.  Corollary 3.7 then applies to yield 

\begin{theorem} Under the above hypothesis,  
the functional
    \[
   {I}(u)=\Psi (u)+\Psi^{*}(-(u\cdot \nabla)u-\vec a\cdot  \nabla u - \frac{1}{2}{\rm div} ( \vec a)u )  
\]
 has zero for infimum and the latter is attained at a solution $\bar u$ for (\ref{NSE5}).  \\
 \end{theorem}     
  We can also give a variational resolution for nonlinear anti-Hamiltonian systems. 

 \begin{theorem} Let $\phi$ be a  proper convex lower semi-continuous function  on $X\times Y$, let $A:X\to Y^{*}$ be any bounded linear operator, let $B_{1}:X\to X^{*}$ (resp., $B_{2}:Y\to Y^{*}$) be two positive bounded linear operators, and assume $\Lambda:=(\Lambda_{1}, \Lambda_{2}): X\times Y \to X^{*}\times Y^{*}$ is a regular conservative operator.  If
 \[
 \lim\limits_{\|x\|+\|y\|\to \infty}\frac{\phi(x,y)+\frac{1}{2} \langle B_{1}x, x\rangle+\frac{1}{2} \langle B_{2}y, y\rangle}{\|x\|+\|y\|}=+\infty,
\]
then for any $(f, g)\in X^{*}\times Y^{*}$, there exists $(\bar x,\bar y) \in X\times Y$ which solves the following system
     \begin{equation}
 \left\{ \begin{array}{lcl}
\label{eqn:existence}
\hfill \Lambda_{1}(x,y)-A^{*} y-B_{1} x +f&\in& \partial_{1} \phi ( x,  y).\\
\hfill \Lambda_{2}(x,y)+A x-B_{2} y+g&\in& \partial_{2} \phi ( x,  y).
\end{array}\right.
 \end{equation}
 
The solution is obtained as a minimizer on $X\times Y$ of the functional 
 \[
I(x,y)=\psi (x, y)+\psi^{*}(-A^{*}y-B^{a}_{1}x+\Lambda_{1}(x,y), Ax-B^{a}_{2}y+\Lambda_{2}(x,y)).
\]
where 
\[
\psi (x,y)=\phi (x,y) +  \frac{1}{2}\langle B_{1}x, x\rangle+\frac{1}{2} \langle B_{2}y, y\rangle -\langle f, x\rangle- \langle g, y\rangle.
\]
and where $B_{1}^{a}$ (resp., $B_{2}^{a}$) are the skew-symmetric parts of $B_{1}$ and $B_{2}$. 
 \end{theorem} 
\noindent {\bf Proof:}  Consider the following ASD Lagrangian  (see \cite{G2})
  \[
L((x,y), (p,q))=\psi(x, y)+\psi^{*}(-A^{*}y-B^{a}_{1}x-p, Ax-B_{2}y-q).
\]
Theorem 4.1 yields that $I(x,y)=L((x,y), \Lambda (x,y))$ attains its minimum at some point $(\bar x,\bar y) \in X\times Y$ and that the minimum is actually $0$. In other words, 
\begin{eqnarray*}
0&=&I(\bar x,\bar y)=\psi(\bar x, \bar y)+\psi^{*}(-A^{*}\bar y-B^{a}_{1}\bar x+\Lambda_{1}(\bar x,\bar y), A\bar x-B^{a}_{2}\bar y+\Lambda_{2}(\bar x,\bar y))\\
&=&\psi (\bar x, \bar y)+\psi^{*}(-A^{*}\bar y-B^{a}_{1}\bar x+\Lambda_{1}(\bar x,\bar y), A\bar x-B^{a}_{2}\bar y+\Lambda_{2}(\bar x,\bar y))\\
&&-\langle (\bar x, \bar y),  (-A^{*}\bar y-B^{a}_{1}\bar x+\Lambda_{1}(\bar x,\bar y), A\bar x-B^{a}_{2}\bar y+\Lambda_{2}(\bar x,\bar y))\rangle
\end{eqnarray*}
from which  follows that
\begin{equation}
 \left\{ \begin{array}{lcl}
\label{eqn:existence}
-A^{*}y-B^{a}_{1}x+\Lambda_{1}(x,y)  &\in& \partial_{1} \varphi (x,y)+ B_{1}^{s}(x)-f\\
\hfill Ax-B^{a}_{2}y+\Lambda_{1}(x,y)+g&\in& \partial_{2} \varphi (x, y) + B_{2}^{s}(y)-g.
\end{array}\right.
\end{equation}

A typical example of such a system are the equations of magneto-hydrodynamics, but here is a simpler example communicated to us by A. Moameni.
\subsection*{Example 4: A variational resolution for a doubly nonlinear coupled equations}
Let ${\bf b_{1}}:\Omega \to  {\bf R^{n}}$ and  ${\bf b_{2}}:\Omega \to  {\bf R^{n}}$ be two smooth vector fields on a bounded domain $\Omega$ of  $\bf R^{n}$, verifying the conditions in example 3 and let $B_{1}v={\bf b_{1}}\cdot \nabla v$ and  $B_{2}v={\bf b_{2}}\cdot \nabla v$ be the corresponding first order linear operators.   Consider the Dirichlet problem:
\begin{equation}
\label{Ex1.500}
 \left\{ \begin{array}{lcl}
    \hfill Ê\Delta (v+u) +{\bf b_{1}}\cdot \nabla u &=& |u|^{p-2}u+u^{m-1}v^{m} +f  \hbox{\rm \, on \, $\Omega$}
\\
  \hfill Ê\Delta (v-c^{2}u) +{\bf b_{2}}\cdot \nabla v&=& |v|^{q-2}q -u^{m}v^{m-1} +g  \hbox{\rm \, on \, $\Omega$}\\
 \hfill  u=v &=& 0 \quad \quad \quad \quad  \hbox{\rm on  $\partial \Omega$. }        \end{array}  \right.
   \end{equation}
 We can use the above to get
 \begin{theorem} Assume  ${\rm div }({\bf b_{1}})\geq 0$ and ${\rm div }({\bf b_{2}})\geq 0$ on $\Omega$, $2<p,q\leq \frac{2n}{n-2}$ and $1<m<\frac{n+2}{n-2}$ and consider on $H^{1}_{0}(\Omega)\times H^{1}_{0}(\Omega)$ 
 the functional 
\[
I(u, v)=\Psi (u) +\Psi^{*}({\bf b_{1}}.\nabla u +\frac{1}{2}{\rm div }({\bf b_{1}}) \,  u+\Delta v-u^{m-1}v^{m} )+\Phi (v) +\Phi^{*}({\bf b_{2}}.\nabla v +\frac{1}{2}{\rm div }({\bf b_{2}}) \,   v -c^{2}\Delta u +u^{m}v^{m-1})
\]
where
\[
 \Psi (u)=\frac{1}{2} \int_{\Omega} | \nabla  u |^{2}dx +\frac{1}{p}\int_{\Omega}|u|^{p} dx+\int_{\Omega}fu dx +\frac{1}{4}\int_{\Omega}{\rm div }({\bf b_{1}}) \, |u|^{2} dx, 
 \]
 \[
 \Phi (v)=\frac{1}{2} \int_{\Omega} | \nabla  v |^{2}dx +\frac{1}{q}\int_{\Omega}|v|^{q} dx+\int_{\Omega}gv dx +\frac{1}{4}\int_{\Omega}{\rm div }({\bf b_{2}}) \, |v|^{2} dx
 \]
and $\Psi^{*}$ and $\Phi^{*}$ are their Legendre transforms. Then there exists $(\bar u, \bar v)\in H^{1}_{0}(\Omega)\times H^{1}_{0}(\Omega)$ such that:
\[
I(\bar u, \bar v)=\inf \{I(u, v);  (u, v)\in H^{1}_{0}(\Omega)\times H^{1}_{0}(\Omega)\}=0, 
\]
and $(\bar u, \bar v)$ is a solution of $(\ref{Ex1.500})$.
 \end{theorem}
\noindent{\bf Proof:} Let $A=\Delta$ on $H^{1}_{0}$, $B_{1}={\rm div }({\bf b_{1}})$, $B_{2}={\rm div }({\bf b_{2}})$ and consider the ASD Lagrangian
\[
L((u, v), (r,s)))=\Psi (u) +\Psi^{*}({\bf b_{1}}.\nabla u +\frac{1}{2}{\rm div }({\bf b_{1}}) \,  u+\Delta v+r )+\Phi (v) +\Phi^{*}({\bf b_{2}}.\nabla v +\frac{1}{2}{\rm div }({\bf b_{2}}) \,   v -c^{2}\Delta u +s).
\] 
It is also easy to verify that the nonlinear operator $\Lambda: H^{1}_{0}\times H^{1}_{0} \to H^{-1}\times H^{-1}$ defined by
\[
\Lambda (u,v)=(-u^{m-1}v^{m}, u^{m}v^{m-1})
\]
is regular and conservative.
 
  \section{Nonlinear evolution equations}
 
 Consider now an evolution triple $X \subset H \subset X^*$, where   $H$ is a Hilbert space with $\braket{}{}$ as scalar product, and where $X$ is a dense vector subspace of $H$, that  is a reflexive Banach space once equipped with its own norm $\| \cdot \|$.   Assuming the canonical injection $X \rightarrow H$, 
continuous, we identify the Hilbert space $H$ with its dual $H^*$ and we
``inject''
$H$ in $X^*$   in such a way that
\[
     \braket{h}{u}_{X^*,X} = \braket{h}{u}_H \quad \hbox{\rm for all $h \in
H$
and all $u \in X$}
\]
This injection  is continuous, one-to-one, and $H$ is also dense in $X^*$. 
 In other words,  the dual $X^*$ of $X$ is represented as the completion of $H$ for the dual norm
$\|h\|=\sup\{{\braket{h}{u}}_H; \|u\|_X \leq 1\}$.

Let $[0,T]$ be a fixed real interval  and   consider the following Banach spaces:
 \begin{itemize} 
\item The space $L^{2}_{X}$ of Bochner integrable functions from $[0,T]$ into $X$ with norm  
\[
     \|u\|^2_{L_2(X)} =(\int_0^T \|u(t)\|_{X}^{2} dt)^{\frac{1}{2}}.
\]
 
\item The space ${\cal X}_2$ of all functions in $L^{2}_{X}$ such that $\dot{u} \in L^{2}_{X^*}$, equipped with the norm 
\[
     \|u\|_{\cal X} = (\|u\|_{L_2(X)}^{2} + \|\dot u\|_{L_2(X^*)}^{2})^{1/2}.
 \]
\end{itemize} 
Note that this last space is different from the Sobolev space
 \[
A^{2}_{X} = \{ u:[0,T] \rightarrow X; \,\dot{u} \in L^{2}_{X}  \}
\]
 and we actually have $A^{2}_{X} \subset {\cal X}_2\subset A^{2}_{X^*}$. 
 
   \begin{definition} A  time dependent Lagrangian on $[0,T]\times X\times X^*$ is any function $L: [0,T]\times X\times X^*\to \R \cup \{+\infty\}$ that is  measurable with respect to   the  $\sigma$-field  generated by the products of Lebesgue sets in  $[0,T]$ and Borel sets in $H\times H$. 
   The Hamiltonian $H_{L}$ of $L$ is the function defined on $[0,T]\times X\times X^*$ by:
   \[
   H(t, x,y)=\sup\{\langle y, p\rangle -L(t, x, p); p\in X^*\}
   \]
  We say that  $L$ is an {\it anti-self dual Lagrangian} (ASD) on $[0,T]\times X\times X^*$  if  for any $t\in [0,T]$, the map $L_t:(x,p)\to L(t, x,p)$ is in $\ASD$: that is if 
  \[
  L^*(t, p,x)=L(t, -x, -p) \quad \hbox{\rm for all $(x,p)\in X\times X^*$}, 
  \]
   where here $L^*$ is the Legendre transform in the last two variables.
  \end{definition}
The most basic time-dependent $ASD$-Lagrangians are again of the form 
\[
L(t,x,p)=\varphi (t, x) +\varphi^{*}(t, -p)
\]
where for each $t$, the function $x\to \varphi (t, x)$ is convex and lower semi-continuous on $X$.  We now show how this property naturally ``lifts'' to path space. For that we associate to each time-dependent Lagrangian $L$ on $[0,T]\times X\times X^*$, the corresponding Lagrangian ${\cal L}$ on the path space $L^{2}_{X}\times L^{2}_{X^*}$ defined by
\[
{\cal L} (u,p):=\int_0^T L(t, u(t), p (t)) dt.
\]
Define the dual of ${\cal L}$ in both variables as
\[
{\cal L}^{*}(q,v) = \sup \left\{ \int_0^T ( \braket{q(t)}{u(t)} + \braket{p(t)}{v(t)} -
     L(t, u(t),p(t))) dt\ ;  (u,p) \in L_X^{2} \times L_{X^*}^{2} \right\}
\]
and denote  the associated Hamiltonian on  path space by:
\[
H_{\cal L}(u,v )=\sup \left\{ \int_0^T ( \braket{p(t)}{v(t)}   -
     L(t, u(t),p(t))) dt\ ;  \, p \in L_{X^*}^{2}  \right\}
\]
The following is standard (see \cite{G2}).

\begin{proposition} Suppose that $L$ is a Lagrangian on $[0,T]\times X\times X^*$, and let ${\cal L}$ be the corresponding Lagrangian on the path space $L^{2}_{X}\times L^{2}_{X^*}$. Then
\begin{enumerate}
\item ${\cal L}^{*}(p,u) = \int_0^T L^{*}(t, p(t),u(t))dt$. 

\item $H_{\cal L}(u,v )=\int_0^T H_{L}(t, u(t),v(t))dt$. 

\item If $L$ is an anti-self dual Lagrangian on $[0,T]\times X\times X^*$, then ${\cal L}$ is anti-selfdual on $L^{2}_{X}$. 

\end{enumerate}
\end{proposition}

  \begin{proposition} 
Suppose $\ell$ is  a self-dual boundary Lagrangian on $H\times H$ and let $L$ be an anti-self dual Lagrangian on $[0,T]\times X\times X^*$ such that  
\begin{equation}
\label{condition.A}
  \hbox{   For each $p\in L^{2}_{X^*}$, the map $u\to \int_0^T L(t, u(t), p (t)) dt$ is continuous on $L^{2}_{X}$ }
 \end{equation}
 \begin{equation}
 \label{condition.B}
  \hbox{   The map $u\to \int_0^T L(t, u(t), 0) dt$ is bounded on the bounded sets of $L^{2}_{X}$ }
 \end{equation}
  \begin{equation}
  \hbox{ $  \ell (a,b) \leq C(1+\|a\|_H^2+\|b\|^2_H)$ for all $(a,b)\in H\times H$.}
   \end{equation}
 Then the Lagrangian
\[
{\cal M}_{L} (u,p)= \left\{ \begin{array}{lcl}
  \hfill  \int_0^T L(t, u(t), p (t)+\dot{u}(t)) dt + \ell (u(0),u(T))\, &{\rm if}& u\in {\cal X}_2 \\
  +\infty &\,& {\rm otherwise} \\
\end{array}\right.
\]
is anti-self dual on $L^{2}_{X}$.
\end{proposition}
 
\noindent{\bf Proof:} 
    For $(q,v) \in L^{2}_{X}\times {\cal X}_2$, write:
 \begin{eqnarray*}
    {\cal M}_{L} ^{*}(q,v) &=& \sup_{u \in L^{2}_{X}} \sup_{p \in
     L^{2}_{X^*}}\{\int_0^T (\braket{u(t)}{q(t)}+\braket{v(t)}{p(t)}
     - L(t, u(t),p(t)+\dot{u}(t) ) )dt -\ell (u(0),u(T)) \}\\
 &=& \sup_{u \in {\cal X}_2} \sup_{p \in
     L^{2}_{X^*}}\{\int_0^T (\braket{u(t)}{q(t)}+\braket{v(t)}{p(t)}
     - L(t, u(t),p(t)+\dot{u}(t) ) )dt -\ell (u(0),u(T)) \}
\end{eqnarray*}
  Make a substitution $
 p(t)+\dot{u}(t)=r(t) \in L_{X^*}^2$.  Since $u$ and $v$ are both in ${\cal X}_2$, we have:
\[
\int_0^T \braket{v}{\dot{u}}=- \int_0^T \braket{\dot v}{u} +
\braket{v(T)}{u(T)} - \braket{v(0)}{u(0)},  \] and since the subspace
  \[
  {\cal X}_{_{2,0}}=\{u\in {\cal X}_2;\,  u(0)=u(T)=0\}
  \]
  is dense in $L^2_X$, we obtain
  \begin{eqnarray*}
  {\cal M}_{L} ^{*}(q,v) &=& \sup_{u \in {\cal X}_2} \sup_{r \in
     L^{2}_{X^*}}\{\int_0^T (\braket{u(t)}{q(t)}+\braket{v(t)}{r(t)-{\dot u}(t)}
     - L(t, u(t),r(t) ) )dt-\ell (u(0),u(T))\\
  &=& \sup_{u \in {\cal X}_2} \sup_{r \in
     L^{2}_{X^*}}\{\int_0^T (\braket{u(t)}{q(t)+{\dot v(t)}}+\braket{v(t)}{r(t)}
     - L(t, u(t),r(t) ) )dt\\
  && \quad \quad \quad \quad \quad \quad \quad \quad-
\braket{v(T)}{u(T)} + \braket{v(0)}{u(0)} - \ell (u(0),u(T))\}\\
&=&  \sup_{u \in {\cal X}_2 } \sup_{r \in
     L^{2}_{X^*}}\sup_{u_0\in {{\cal X}_{2,0}}}\{\int_0^T (\braket{u(t)}{q(t)+{\dot v(t)}}+\braket{v(t)}{r(t)}
     - L(t, u(t),r(t) ) )dt\\
  && \quad \quad \quad \quad \quad \quad \quad \quad-
\braket{v(T)}{(u+u_0)(T)} + \braket{v(0)}{(u+u_0)(0)} - \ell ((u+u_0)(0),(u+u_0)(T))\}\\
  &=&  \sup_{w \in {\cal X}_2 } \sup_{r \in
     L^{2}_{X^*}}\sup_{u_0\in {{\cal X}_{2,0}}}\{\int_0^T (\braket{w(t)-u_0(t)}{q(t)+{\dot v(t)}}+\braket{v(t)}{r(t)}
     - L(t, w(t)-u_0(t),r(t) ) )dt\\
  && \quad \quad \quad \quad \quad \quad \quad \quad-
\braket{v(T)}{w(T)} + \braket{v(0)}{w(0)} - \ell (w(0), w(T))\}\\
  &=&  \sup_{w\in {\cal X}_2 } \sup_{r \in
     L^{2}_{X^*}}\sup_{x\in L^2_X}\{\int_0^T (\braket{x(t)}{q(t)+{\dot v(t)}}+\braket{v(t)}{r(t)}
     - L(t, x(t),r(t) ) )dt\\
  && \quad \quad \quad \quad \quad \quad \quad \quad -
\braket{v(T)}{w(T)} + \braket{v(0)}{w(0)} - \ell (w(0), w(T))\}
  \end{eqnarray*}
Here we have used the fact  that ${\cal X}_{2,0}$ is dense in $L^{2}_{X}$ and the continuity of $u\to \int_{0}^{T}L(t, u(t), p(t)) dt$ on $L^{2}_{X}$ for each $p$. 

Now, for each $(a,b)\in X\times X$, there is $w\in {\cal X}_2$ such that $w(0)=a$ and $w(T)=b$, namely the linear path $w(t)=\frac{(T-t)}{T}a +\frac{t}{T}b$. Since also $X$ is dense in $H$ and $\ell$ is continuous on $H$,  we finally obtain that 
  \begin{eqnarray*}
 {\cal M}_{L} ^{*}(q,v) &=&\sup_{(a,b)\in X\times X} \sup_{r \in
     L^{2}_{X^*}}\sup_{x\in L^2_X}\{\int_0^T (\braket{x(t)}{q(t)+{\dot v(t)}}+\braket{v(t)}{r(t)}
     - L(t, x(t),r(t) ) )dt\\
  && \quad \quad \quad \quad \quad\quad\quad\quad\quad \quad\quad\quad\quad -
\braket{v(T)}{b} + \braket{v(0)}{a} - \ell (a, b)\}\\
  &=& \sup_{x \in L^{2}_{X}} \sup_{r \in
     L^{2}_{X^*}}\{\int_0^T (\braket{x(t)}{q(t)+{\dot v(t)}}+\braket{v(t)}{r(t)}
     - L(t, x(t),r(t) ) )dt\\
   && +
  \sup_{a\in H}\sup_{b\in H}\{-\braket{v(T)}{b} + \braket{v(0)}{a} - \ell (a,b)\} \\
  &=& \int_0^T L^{*}(t, q(t)+\dot{v}(t), v(t)) dt + \ell^{*} (v(0),-v(T))\\
  &=& \int_0^T L(t, -v(t), -\dot{v}(t)-q(t)) dt + \ell (-v(0),-v(T))\\
  &=& M(-v,-q). 
\end{eqnarray*}
If now $(q,v) \in L^{2}_{X^*}\times L^{2}_{X}\setminus {\cal X}_2$, then we use the fact that $u\to \int_{0}^{T}L(t, u(t), 0) dt$ is bounded on the unit ball of ${\cal X}_2$ and the growth condition on $\ell$ to deduce
 \begin{eqnarray*}
    {\cal M}_{L} ^{*}(q,v) &\geq & \sup_{u \in {\cal X}_2} \sup_{r \in
     {\cal X}_2}\{\int_0^T (\braket{u(t)}{q(t)}+\braket{v(t)}{r(t)}-\braket{v(t)}{\dot{u}(t)}
     - L(t, u(t), r(t) ) )dt -\ell (u(0),u(T)) \}\\
 &\geq&  
 \sup_{u \in {\cal X}_2} \sup_{r \in
     {\cal X}_2}\{-\|u\|_{L^2_X}\|q\|_{L^2_{X^*}}-\|v\|_{L^2_X}\|r\|_{L^2_{X^*}}+\int_0^T (-\braket{v(t)}{\dot{u}(t)}
     - L(t, u(t), r(t) ) )dt -\ell (u(0),u(T)) \}\\
 &\geq&  
 \sup_{\|u\|_{{\cal X}_2}\leq 1} \{-\|q\|_{2} +\int_0^T (\braket{-v(t)}{\dot{u}(t)}
     - L(t, u(t), 0) )dt -\ell (u(0),u(T)) \}\\
 &\geq&  
 \sup_{\|u\|_{{\cal X}_2}\leq 1}  \{C+\int_0^T (\braket{-v(t)}{\dot{u}(t)}
     - L(t, u(t), 0 ) )dt - \frac{1}{2}(\|u(0)\|^2 +\|u(T)\|^2)\}\\
  &\geq&  
 \sup_{\|u\|_{{\cal X}_2}\leq 1}  \{D+\int_0^T \braket{-v(t)}{\dot{u}(t)}dt - \frac{1}{2}(\|u(0)\|_X^2 +\|u(T)\|_X^2)\}.
\end{eqnarray*}
Since now $v$ does not belong to ${\cal X}_2$, we have that 
\[
 \sup_{\|u\|_{{\cal X}_2}\leq 1}\int_0^T (\braket{v(t)}{\dot{u}(t)}dt + \frac{1}{2}(\|u(0)\|_X^2 +\|u(T)\|_X^2)=+\infty 
 \]
 which means that $M^{*}(q,v)=+\infty = M(-v, -q)$. \\

Now we can prove the following
 \begin{theorem} Let $X\subset H \subset X^*$ be an evolution pair and consider an anti-self dual Lagrangian   $L$   on $[0,T]\times X\times X^*$ and a self-dual boundary Lagrangian $\ell $ on $H\times H$. Assume the following conditions:
    \begin{equation} 
    \label{Condition.A1}
 \hbox{   For each $p\in L^{2}_{X^*}$, the map $u\to \int_0^T L(t, u(t), p (t)) dt$ is bounded on the bounded sets of $L^{2}_{X}$ }
 \end{equation}
  \begin{equation}
   \label{Condition.A2}
     \lim\limits_{\|v\|_{L^2(X)}\to +\infty}\int_{0}^{T}H_{L}(t, 0, v(t)) dt=+\infty, 
   \end{equation}
   and
    \begin{equation}
     \label{Condition.A3}
  \hbox{ $  \ell (a,b) \leq C(1+\|a\|_H^2+\|b\|^2_H)$ for all $(a,b)\in H\times H$.}
   \end{equation}
(1) Then for any regular conservative operator $\Lambda: D(\Lambda)\subset L^2(X)\to L^2(X^*)$ such that ${\cal X}_2\subset D(\Lambda)$,  the following functional
    \[
    I_{\ell,L, \Lambda}(u) = \int_0^T L(t, u(t),\Lambda u(t)+\dot{u}(t))dt + \ell (u(0),u(T))
    \]
   has zero infimum. Moreover,  there exists $v \in {\cal X}_2$ such that:
   \begin{equation}
   \label{eqn:2.18}
\hbox{\rm  $\big(v(t),\Lambda v(t)+\dot{{v}} (t)\big)\in \mbox{\rm Dom} (L)$ for almost all 
$t\in [0,T]$}
\end{equation} 
\begin{equation}
\label{eqn:2.19}
I_{\ell,L, \Lambda}( v)=\inf\limits_{u\in {\cal X}_2}I_{\ell, L, \Lambda}(u)=0,
\end{equation}
\begin{equation}
\label{eqn:2.20}
\hbox{\rm $L(t, v(t),\Lambda v(t)+\dot{v}(t))+\langle v(t),\dot v(t)\rangle =0$ for almost all $t\in [0,T]$},
\end{equation}
\begin{equation}
\label{eqn:2.21}
\ell (v(0), v(T))=\frac{1}{2}(\|v(T)\|^2_H-\|v(0)\|^2_H),  
\end{equation}
 \begin{equation}
\label{eqn:2.22}
 (-{\dot v}(t)-\Lambda v(t), -v(t))\in \partial L(t, v(t), {\dot v} (t)+\Lambda v(t)).
\end{equation}
 (2)   In particular, for every $v_{0}\in H$ the following functional 
\[
 I_{v_{0},L,\Lambda}(u)=  \int_0^T L(t, u(t),\Lambda u(t)+\dot{u}(t))dt +\frac{1}{2}\|u(0)\|^{2} -2\langle v_0, u(0)\rangle +\|v_0\|^{2} +\frac{1}{2}\|u(T)\|^{2}
\]
has minimum equal to zero on $L^2_X$. It  is attained at a unique path $v$ such that $v(0)=v_{0}$, verifying  (\ref{eqn:2.18}- \ref{eqn:2.22})
and in particular
\begin{equation}
\label{eqn:2.23}
\hbox{\rm $ \|v(t)\|_{H}^{2}=\|v_{0}\|^{2}-2\int_{0}^{t}L(s, v(s), \Lambda v(s)+{\dot v}(s)) ds$   for every $t\in [0,T]$.}
\end{equation}
 \end{theorem} 
 \noindent{\bf Proof:} First apply Proposition 5.2 to get that the Lagrangian
\[
{\cal M}_{L} (u,p)= \left\{ \begin{array}{lcl}
  \hfill  \int_0^T L(t, u(t), p (t)+\dot{u}(t)) dt + \ell (u(0),u(T))\, &{\rm if}& u\in {\cal X}_2 \\
  +\infty &\,& {\rm otherwise} \\
\end{array}\right.
\]
is anti-self dual on $L^{2}_{X}$.  It is now sufficient to apply Corollary 3.5 to conclude that the infimum of ${\cal M}_{L} (u,\tilde \Lambda u)$ is equal $0$ and is achieved. This yields claim (\ref{eqn:2.18}) and (\ref{eqn:2.19}).

 Since $ L(t,v(t),\dot v(t))\geq - \langle v(t), \dot v(t) \rangle $ for all $t\in [0,T]$,  and since  $\ell (v(0), v(T))\geq \frac{1}{2}(\|v(T)\|^2_H-\|v(0)\|^2_H)$, claims (\ref{eqn:2.20}) and (\ref{eqn:2.21}) follow from the following identity
 \[
 0=   I_{\ell,L, \Lambda}(v) = \int_0^T L(t, v(t),\Lambda v(t)+\dot{v}(t)+ \langle v(t),\dot v(t)\rangle)dt -\frac{1}{2}(\|v(T)\|^2_H-\|v(0)\|^2_H)+ \ell (v(0),v(T)).
 \]
  To prove (\ref{eqn:2.22}), use (\ref{eqn:2.20}), the fact that $L$ is anti-selfdual and that $\Lambda$ is conservative to write:
 \[
  L(s, v(s),\Lambda v(s)+\dot{v}(s))+L^{*}(s, -\Lambda v(s)-\dot{v}(s), -v(s)) +\langle (v(s), \Lambda v(s)+{\dot v}(s)), (\Lambda v(s)+{\dot v}(s), v(s) )\rangle=0
 \]   
and conclude by the limiting case of the Legendre-Fenchel  duality in the space $X\times X^*$.  \\ 
For (2) it suffices to apply the first part with   the boundary Lagrangian 
\[
\ell (r,s)= \frac{1}{2}\|r\|^{2} -2\langle v_{0}, r\rangle +\|v_{0}\|^{2} +\frac{1}{2}\|s\|^{2}. 
\]
which is clearly self-dual. We then get
 \[
I_{\ell,L,\Lambda}(u) = \int_0^T \left[ L(t, u(t),\Lambda u(t)+\dot{u}(t))+\langle u(t), {\dot u}(t)\rangle \right] dt +  \|u(0)-v_{0}\|^{2}.
\]
Note also that (\ref{eqn:2.20}) yields
   \[
   \frac{d (|v(s)|^{2})}{ds}=-2L(s, v(s),\Lambda v(s)+{\dot v}(s)). 
   \]
which readily implies  (\ref{eqn:2.23}).\\

 We now apply the results of the last section to the particular class of ASD Lagrangian of the form $L(x,p)=\phi (x)+\phi^*(Ax-p)$ to obtain variational formulations and proofs of existence for various nonlinear  parabolic equations.

 \begin{proposition}
\label{Asym flow}
Let $X\subset H \subset X^*$ be an evolution triple and consider for each $t\in [0,T]$ a bounded linear operator $A_t:X\to X^*$ and  $\phi : [0,T]\times X\to \bar\R$ such that for each $t$ the functional $\psi (t, x):=\phi (t, x) +\frac{1}{2} \langle A_tx, x\rangle$ is  convex, lower semi-continuous and satisfies for some $C>0$, $m,n>1$ the following growth condition:
\begin{equation}
\label{growth}
\hbox{$\frac{1}{C}\left( {\| x\|}_{L^2_X}^m-1\right)\le\int_0^T\{\phi (t, x(t)) +\frac{1}{2} \langle A_tx(t), x(t)\rangle\}\, dt\le C\left( {\| x\|}_{L^2_X}^n+1\right)$ for every $x\in L^2_X$.}
\end{equation}
If $\Lambda: D(\Lambda)\subset X \to X^*$ is a regular conservative operator and $v_0\in X$, we consider on  ${\cal X}_2$ the functional
 \[
 I(x) = \int_0^T\{\psi(t, x(t))+\psi^*(t, -\Lambda x(t)-A_t^a x(t) - \dot{x}(t))\}\, dt  + \frac{1}{2}(|x(0)|^{2}+|x(T)|^{2})
 -2\langle x(0), v_{0} \rangle +|v_{0}|^{2}, 
\]
 where for each $t\in [0,T]$, $A^a_t$ is the anti-symmetric part of the operator $A_t$. 
Then there exists a path $v \in {\cal X}_2$ such that
\begin{equation}
 I(v)=\inf\limits_{x\in {\cal X}_2 }I(x)=0.
\end{equation}
 \begin{eqnarray}
 \label{eqn:1000}
-\dot{v}(t) - A_t v(t) -\Lambda  v(t) &\in&\partial\phi(t, v(t))\quad \mbox{ for a.e. }t\in [0,T] \\
 v(0)& = &v_0. \nonumber
 \end{eqnarray}
 \end{proposition}
\noindent{\bf Proof:} The Lagrangian $L(t, x,p) :=\psi(t, x) +\psi^*(t, -A^ax - p)$  is an ASD Lagrangian on $X\times X^*$ by Proposition 2.5. Consider $\ell$ on $H\times H$ to  be 
$
\ell (r,s)=\frac{1}{2}(|r|^{2}+|s|^{2})-2\langle r, v_{0} \rangle +|v_{0}|^{2}, 
$
and lift $\Lambda$ to a regular conservative operator $\tilde \Lambda$ from its domain in $L^2_X([0,T])$ into $L^2_{X^*}([0,T])$ by setting $(\tilde \Lambda x)(t)=\Lambda (x(t))$.
 It is easy to check that all the conditions of Theorem 5.2 are satisfied by $L$, $\ell$, $B$ and $\tilde \Lambda$, hence  there exists $ v\in {\cal X}_2 $ such that $I(v)=0$. We obtain
 
\[
 0=\int_0^T\big(\psi(t, v(t))+\psi^*(t, -\Lambda v(t)-A_t^a v(t) - \dot{v}(t))+ \langle v(t), \Lambda v(t)+A_t v(t) + \dot{v}(t)\rangle\big)\, dt
  +  \frac{1}{2}\Vert v(0)-v_0\Vert_H^2
\]
which yields since  the integrand is non-negative for each $t$ and since we are now in the limiting case of Legendre-Fenchel duality that
 \begin{eqnarray}
 \label{eqn:1000}
-\dot{v}(t) - A^a_t v(t) -\Lambda  v(t) &\in&\partial\phi(t, v(t))+A^s_t v(t)\quad \mbox{ for a.e. }t\in [0,T] \\
 v(0)& = &v_0. \nonumber
\end{eqnarray}

  \subsection*{Example 5: Navier-Stokes evolutions} 
We now consider the evolution equation associated to a fluid driven by its boundary.
\begin{eqnarray}
\label{NSE50}
 \left\{ \begin{array}{lcl}
    \hfill
\frac{\partial u}{\partial t}+(u\cdot \nabla)u +f &=&\nu \Delta u - \nabla  p \quad \hbox{\rm on $[0,T]\times \Omega$}\\
\hfill {\rm div} u&=&0 \quad \quad \quad \hbox{\rm on $[0,T]\times \Omega$}\\
\hfill u(t,x)&=&u^0(x) \quad \hbox{\rm on $[0,T]\times \partial \Omega$}\\
\hfill u(0,x)&=&u_{0}(x) \quad \hbox{\rm on $\Omega$} 
\end{array}\right.
\end{eqnarray}
 where $\int_{\partial \Omega} u^{0} {\bf \cdot n}\, d\sigma =0$, $\nu>0$ and $f\in L^{p}(\Omega;\R^{3})$.  Assuming  that $u^{0} \in H^{3/2}(\partial \Omega)$ and that $\partial \Omega$ is connected,  Hopf's extension theorem again yields  the existence of $v^0\in H^{2}(\Omega)$ such that
\begin{equation}
\label{Hopf}
v^0=u^{0}\,\, \hbox{\rm on $\partial \Omega$,\quad  ${\rm div}\, v^0=0$\quad  and \quad  $\int_{\Omega}\Sigma_{j,k=1}^{3}u_{k}\frac {\partial v^0_{j} } {\partial x_{k}}u_{j}\, dx\leq \epsilon \|u\|^{2}_{V}$ for all $u\in V$}
\end{equation}
where  $V=\{u\in H^{1}(\Omega; {\bf R}^{3}); {\rm div} v=0\}$.
Setting $v=u+v^0$, then solving  (\ref{NSE50}) reduces to finding a solution in the Banach space $V_0=\{u\in H^{1}_0(\Omega; {\bf R}^{3}); {\rm div} v=0\}$
 for
 \begin{eqnarray}
\label{NSE51}
\frac{\partial u}{\partial t}+ (u\cdot \nabla)u +(v^0\cdot \nabla)u +(u\cdot \nabla)v^0 +g   &\in& -\partial \Phi (u) \\
\hfill  u(0)&=& u_0-v^0. 
 \nonumber 
\end{eqnarray}
where $\Phi$ is again the convex Dirichlet energy functional $\Phi (u)=\frac{\nu}{2} \int_{\Omega}\Sigma_{j,k=1}^{3}(\frac {\partial u_{j} } {\partial x_{k}})^{2}\, dx$  and where 
\[
g:= f -\nu \Delta v^0 + (v^0\cdot \nabla)v^0\in V^{*}.
\]
In other words, this is an equation of the form
\begin{equation}
 \label{Lax.Mil.again}
\frac{\partial u}{\partial t}+  \Lambda u+Bu+g\in -\partial \Phi (u) 
  \end{equation}
  where $\Lambda u= (u\cdot \nabla)u$ is a regular conservative operator, and $Bu=(v^0\cdot \nabla)u +(u\cdot \nabla)v^0$ is a bounded linear operator on $V$. The component $B^1u:=(v^0\cdot \nabla)u$ of $B$ is skew-symmetric  which means that  Hopf's estimate implies  
  \[
C\|u\|_V^2\geq  \Psi (u):=\Phi (u)+\frac{1}{2}\langle Bu,u\rangle \geq \frac{1}{2}(\nu-\epsilon)\|u\|^{2} \quad {\rm for\,  all} \,\, u\in V.
  \]
Letting $A^a$ be the antisymmetric part of the operator $Au=(u\cdot \nabla)v^0$, we  can now apply Proposition 5.3  to obtain
  \begin{theorem} Under the above hypothesis on $u^0$,  and for $f\in L^{p}(\Omega, \R^{3})$ with $p>\frac{6}{5}$ and $u_0\in V$,  the minimum of the functional
 
\begin{eqnarray*}
{I}(u)&=&\int_{0}^{T}\left\{\Psi (u(t))+\Psi^{*}(-(u\cdot \nabla)u-B^au+f -\dot u)-\int_{\Omega}\langle f, u\rangle dx\right\} dt\\
&& \quad \quad + \int_{\Omega}\left\{\frac{1}{2}(|u(0,x)|^{2}+|u(x, T)|^{2})-2\langle u(0,x), u_0(x)-v^0(x)\rangle +|u_0(x)-v^0(x)|^{2}\right\}dx
\end{eqnarray*}
 on $A^{2}_{V}$ is zero and is attained at a solution of the equation (\ref{NSE51}). 
\end{theorem}
  
\section{Autonomous ASD Lagrangians and nonlinear conservative operators}

In order to deal with evolutions involving unbounded operators, one has to relax the stringent boundedness condition (\ref{Condition.A1}) on the Lagrangian $L$ since it must be assigned infinite values in some part of the space. In order to do that, we shall $\lambda$-regularize $L$ via inf-convolution to obtain an ASD Lagrangian $L_\lambda$ to which Theorem 5.3 applies and then we let $\lambda \to 0$ in a way reminescent of  Yosida's regularization of unbounded operators and convex functions. 

We shall be able to carry this program in the case where the Lagrangian $L(x, p)$ is autonomous and we obtain the following result. First we  define the {\it Partial Domain} of $\partial L$ to be the set:
\[
{\rm Dom}_{1}(\partial L)=\{x\in H;\hbox{\rm  there exists $p,q\in H$ such that $(p,0)\in \partial L(x, q)$\}.}
\]
 Note that if $L(x,p)=\phi (x) +\phi^{*}(-p)$ with $0$ assumed to be in the domain of $\partial \phi$, then $x_0$ belongs to  ${\rm Dom}_{1}(\partial L)$ if and only if it belongs to the domain of $\partial \phi$. 
This is the main result of this section 

\begin{theorem} Let $X\subset H \subset X^*$ be an evolution triple and let $L$ be an  anti-selfdual Lagrangian on $X\times X^*$ that is uniformly convex in the first variable on $H$ and such that 
 \begin{equation} 
      \label{Compa.1}
    \hbox{For all $x\in X$, the map $L(x,\cdot ):X^*\to\overline{\R}$ is  continuous on $X^*$.}
     \end{equation}
    \begin{equation}
       \label{Compa.2}
    \hbox{ There exists $x_0\in X$ such that $p\to L(x_0,p)$ is bounded on the bounded sets of $X^*$.}
    \end{equation}

Assuming ${\rm Dom}_{1}(\partial L)$ is non-empty, then for every $x_0\in {\rm Dom}_{1}(\partial L)$ there exists  a path $x(t)\in A^2_H$ satisfies $x(0)=x_0$ and 
 \begin{equation}
\label{eqn:2.2}
 (-\Lambda x(t)-{\dot x}(t), -x(t))\in \partial L(x(t), \Lambda x(t)+{\dot x} (t))
\end{equation}
  It is obtained as a minimizer on $A_{H}^2$ of the functional
\[
 I(u)=  \int_0^T L(u(t),\Lambda u(t)+\dot{u}(t))dt +\frac{1}{2}\|u(0)\|^{2} -2\langle x, u(0)\rangle +\|x\|^{2} +\frac{1}{2}\|u(T)\|^{2}.
\]
 \end{theorem}
 We need the following notions. \\

\begin{definition} Let $X\subset H \subset X^*$ be  an evolution triple and let $L$ be an ASD  Lagrangian  on $X\times X^*$. We say that 
\begin{enumerate}
\item   $L$ is $H$-compatible if it lifts  to an ASD Lagrangian on  $H\times H$. In other words if  the Lagrangian 
\begin{eqnarray*}
\tilde L(x,p):=\left\{ \begin{array}{ll}L(x,p) &x\in X\\
+\infty &x\in H\backslash X\end{array}\right.
\end{eqnarray*}
is also anti selfdual on $H\times H$.
 
\item  $L$  is {\it uniformly convex in the first variable (resp. second variable)} on $H$  if there exists $\var >0$ such that for all $p\in H$ (resp. for all $x\in H$) the Lagrangian 
\[
 {\tilde L}(x,p)-\frac{\var\| x\|_H^2}{2} \left(\mbox{resp. }L(x,p)-\frac{\var\| p\|_H^2}{2}\right)
 \]
 is convex in $x$ (resp. in $p$) on $H$.   
 \end{enumerate}
\end{definition}

We start by proving the following proposition which improves on Theorem 5.3 in the case of autonomous Lagrangians.
\begin{proposition}
\label{lift}
 Let $X\subset H \subset X^*$ be an evolution pair and consider a self-dual boundary Lagrangian $\ell $ on $H\times H$. Suppose $L$ is an $H$-compatible autonomous anti-self dual Lagrangian   $L$   on $ X\times X^*$ that is uniformly convex in both variables and satisfying the following conditions:
    \begin{equation} 
    \label{Cond.1}
 \hbox{   For each $p\in L^{2}_{X^*}$, the map $u\to \int_0^T L(u(t), p (t)) dt$ is bounded on the bounded sets of $L^{2}_{X}$ }
 \end{equation}
  \begin{equation}
    \label{Cond.2}
     \lim\limits_{\|v\|_{L^2(X)}\to +\infty}\int_{0}^{T}H_{L}(0, v(t)) dt=+\infty, 
   \end{equation}
    \begin{equation}
      \label{Cond.3}
  \hbox{ $  \ell (a,b) \leq C(1+\|a\|_H^2+\|b\|^2_H)$ for all $(a,b)\in H\times H$.}
   \end{equation}
 
  Then for any regular conservative operator $\Lambda: D(\Lambda)\subset L^2(X)\to L^2(X^*)$ such that ${\cal X}_2\subset D(\Lambda)$,  the following functional
    \[
    I_{\ell,L, \Lambda}(u) = \int_0^T L(t, u(t),\Lambda u(t)+\dot{u}(t))dt + \ell (u(0),u(T))
    \]
   has zero infimum. Moreover,  there exists $v \in {\cal X}_2\cap C^1([0,T], H)$ such that:
   \begin{equation}
   \label{eqn:20.18}
\hbox{\rm  $\big(v(t),\Lambda v(t)+\dot{{v}} (t)\big)\in \mbox{\rm Dom} (L)$ for   all 
$t\in [0,T]$}
\end{equation} 
\begin{equation}
\label{eqn:20.19}
I_{\ell,L, \Lambda}( v)=\inf\limits_{u\in {\cal X}_2}I_{\ell, L, \Lambda}(u)=0,
\end{equation}
\begin{equation}
\label{eqn:20.20}
\hbox{\rm $L(t, v(t),\Lambda v(t)+\dot{v}(t))+\langle v(t),\dot v(t)\rangle =0$ for   all $t\in [0,T]$},
\end{equation}
\begin{equation}
\label{eqn:20.21}
\ell (v(0), v(T))=\frac{1}{2}(\|v(T)\|^2_H-\|v(0)\|^2_H),  
\end{equation}
 \begin{equation}
\label{eqn:20.22}
 (-{\dot v}(t)-\Lambda v(t), -v(t))\in \partial L(t, v(t), {\dot v} (t)+\Lambda v(t)).
\end{equation}
There is $C>0$ such that 
\begin{equation}
\label{time.estimate}
\|\dot v(t)\|\leq C\|\dot v(0)\| \quad \hbox{\rm  for all $t\in [0,T]$.}
\end{equation}
  \end{proposition} 
\noindent {\bf Proof:} Note that Theorem 5.2 already gives the existence of $v\in {\cal X}_2$ satisfying (\ref{eqn:20.18})- (\ref{eqn:20.22}). We shall use   the fact that $L$ is $H$-compatible and uniformly convex in both variables to get the last estimate (\ref{time.estimate}). 
  
Indeed,  since $\tilde L$  is an ASD Lagrangian on $H\times H$ that is uniformly convex in both variables, we can apply Lemma 4.2 of \cite{GT2} to conclude that that $(x,p)\to \partial { \tilde L}(x,p)$ is Lipschitz  on $H\times H$ and Lemma 4.4 of \cite{GT2} to get that $v\in C^1([0,T], H)$.  It follows by continuity that  
\begin{equation}
   \label{Cond.7}
(-{\dot v}(t)-\Lambda v(t), -v(t))\in \partial L(t, v(t), {\dot v} (t)+\Lambda v(t))
\end{equation}
 holds for all $t\in [0,T]$. 

To establish (\ref{time.estimate}),  we first differentiate to obtain:
\[
\frac{1}{2}\frac{d}{dt} \|v(t)-v(t+h)\|^2=
\braket{v(t)-v(t+h)}{\dot v(t)-\dot v(t+h)}.
\]
Setting 
$v_1(t) =\partial_1L\big( v(t), \Lambda v(t) + \dot v(t)\big)$ and 
$v_2(t) =\partial_2L\big(v(t), \Lambda v(t)+\dot v(t)\big)$, we obtain 
 from (\ref{eqn:20.22})  and monotonicity that 
\begin{eqnarray*}
\frac{d}{dt} \| v(t)- v(t+h)\|^2&=&\langle v(t)-v(t+h),-v_1(t)-\Lambda v (t)+v_1(t+h)+\Lambda v (t+h)\rangle\\
&& +
 \langle \dot v(t)+\Lambda v(t) -\Lambda v(t)-\dot v(t+h)-\Lambda v(t+h) +\Lambda v(t+h),-v_2(t)+v_2(t+h)\rangle\\
&\leq &\langle v(t+h),\Lambda v (t)\rangle +\langle v(t), \Lambda v (t+h)\rangle  +\langle -\Lambda v(t)+\Lambda v(t+h), -v_2(t)+v_2(t+h)\rangle.
\end{eqnarray*}
Since $v\in C^1([0,T], H)\cap {\cal X}_2$, we have that for each $t\in [0,T]$, $\lim_{t\to 0} v(t+h)-v(t)=0$ weakly in $X$.  Since $\Lambda$ is regular conservative, it follows that 
\[
\limsup_{t\to 0}\frac{d}{dt} \| v(t)- v(t+h)\|^2 \leq 0.
\]
from which (\ref{time.estimate}) follows. \\
 
The strategy for the proof of Theorem 6.1, is to first consider the $\la$-regularized Lagrangian $L_\la$ that will satisfy teh conditions of Proposition 6.1 and then to try to conclude by letting $\la$ go to zero.  We shall first summarize the needed  properties about inf-convolution of ASD Lagrangians, many of which were established in details in \cite{GT2}.  Suppose $L$ is a Lagrangian on $H\times H^*$ and 
 recall that for each $\la >0$, 
\begin{eqnarray*}
L_\la (x,p):=(L\star T_\lambda)(x,p)=\inf_z\left\{ L(z,p)+\frac{\| x-z\|^2}{2\la}\right\} +\frac{\la\| p\|^2}{2}
\end{eqnarray*}
We have seen (Proposition ??) that if $L$ is an ASD Lagrangian,  then $L_{\lambda}$ is then  a tempered  ASD Lagrangian.  We let now $J_\la (x,p)$ be the unique minimizer in $H$ of the following optimization problem
\[
 \inf_z\left\{ L(z,p)+\frac{\| x-z\|^2}{2\la}\right\}
 \]
 in such a way that
 \begin{equation}
 \label{lambda.formula}
L_\la (x,p)= L(J_\la (x,p),p)+\frac{\| x-J_\la (x,p)\|^2}{2\la} +\frac{\la\| p\|^2}{2}. 
\end{equation}
It is clear that 
\begin{equation}
\label{lambda.derivative}
 \partial_1L_\la (x,p)=\frac{x-J_\la (x,p)}{\la}\in \partial_1L\big( J_\la (x,p),p\big).
 \end{equation}
The following proposition summarizes various properties of this regularization procedure on ASD Lagrangians. For the proofs we refer to \cite{GT2}.  
\begin{lemma} Let $L$ be an ASD Lagrangian on a Hilbert space $H\times H$. 
\begin{enumerate}
\item  If  $L$ is uniformly convex  in the first variable on $H$, then for each $ \la >0$, the Lagrangian $L_\la$ is uniformly convex in both variables on $H\times H$.
  \item If  $L$ is an ASD Lagrangian that is uniformly convex in the first variable, then the map $(x,p)\to J_\la (x,p)$ is Lipschitz from $H\times H$ into $H$.
  
 \item If $L$ is an  ASD Lagrangian on a Hilbert space $H$, then if $(x,p)$ satisfy  $-(p,x)=\partial L_\la (x,p)$ then  
$-\big( p, J_\la (x,p)\big)\in \partial L\big( J_\la (x,p),p\big)$.

\item If $0\in {\rm dom}_{1}\partial L$ where $L$ is an ASD Lagrangian on a Hilbert space, then there exists a constant $C>0$, such that whenever 
 $y_\la$ satisfies $ -(y_\la ,0)\in \partial L_\la (0,y_\la )$ then $\| y_\la\|\le C$ for all $\la >0$.
\end{enumerate}
\end{lemma}

\begin{lemma} Let $X\subset H \subset X^*$ be an evolution triple and let $L$ be an ASD Lagrangian on $X$ that verifies  conditions (\ref{Compa.1}) and (\ref{Compa.2}). 
 Then $L$ is an $H$-compatible ASD Lagrangian. 
\end{lemma}

\noindent {\bf Proof:} We check that $\tilde L$ is still anti-selfdual on $H\times H$.  For $(\tilde x,\tilde p)\in X\times H$, write
\begin{eqnarray*}
{\tilde L}^*(-\tilde p,-\tilde x) &=& \sup_{\stackrel{x\in X}{p\in H}}\left\{ {\braket{\tilde x}{p}}_H 
  +{\braket{\tilde p}{x}}_H-L(x,p)\right\}\\
&=&\sup_{x\in X}\sup_{p\in H}\left\{{\braket{\tilde x}{p}}_{X,X^*}
  + \braket{x}{\tilde p}_{X, X^*}
   -L(x,p)\right\}\\
&=&\sup_{x\in X}\sup_{p\in X^*}\left\{\braket{\tilde x}{p}_{X, X^*}
  +\braket{x}{\tilde p}_{X,X^*}
   -L(x,p)\right\}\\
 &=& L(-\tilde x,-\tilde p)
\end{eqnarray*}
Now suppose $\tilde x\in H\backslash X$. Then
 \[
{\tilde L}^*(\tilde p,\tilde x) = \sup_{\stackrel{x\in X}{p\in H}} \left\{ {\braket{\tilde x}{p}}_H
   +{\braket{\tilde p}{x}}_H
   -L(x,p)\right\}\\
\ge \braket{\tilde p}{x_0} +\sup_{p\in H}\left\{ {\braket{\tilde x }{p}}_H
   -L(x_0,p)\right\}. 
\]
Since $\tilde x\notin X$,  we have that $\sup \left\{ \braket{\tilde x }{p}; p\in H, \|p\|_{X^*}\leq 1\right\} =+\infty $. Since $p\to L(x_0,p)$ is bounded on the bounded sets of $X^*$, it follows that 
$
{\tilde L}^*(\tilde p,\tilde x) \geq  \braket{\tilde p}{x_0} +\sup_{p\in H}\left\{ {\braket{\tilde x }{p}}_H
   -L(x_0,p)\right\}=+\infty,
$
and therefore ${\tilde L}$ is an ASD Lagrangian on $H$.\\

\noindent{\bf Proof of Theorem 6.1:} Since $L$ satisfies (\ref{Compa.1}) and (\ref{Compa.2}), we get from the preceeding lemma that it lifts to an ASD-Lagrangian $\tilde L$ on $H\times H$. For each  $\la >0$, we denote by $L_\la$ the $\lambda$-regularization of $\tilde L$. It  satisfies all the hypothesis of Proposition 6.1, hence 
there exists  then $x_\la\in C_H^1\big([0,T]; H)$ such that  $x_\la (0)=x_{0}$
\begin{equation}
 \label{sol1}
\int_0^T L_\la \big(x_\la (t),\Lambda x_{\lambda}(t)+\dot x_\la (t)\big)\, dt+\ell\big( x_\la (0),x_\la (T)\big) =0
 \end{equation}
\begin{equation}
 \label{sol2}
 \big(-\Lambda x_{\lambda}(t)-\dot x_\la (t),x_\la (t)\big)\in \partial L_\la \big(x_\la (t),\Lambda x_{\lambda}+\dot x_\la (t)\big) 
 \end{equation}
\begin{equation}
 \label{sol3}
 \|\dot x_\la (t)\|\le C(T)\|\dot x_\la (0)\|\quad \hbox{\rm  for all $t\in [0,T]$}.
\end{equation}
Now recall that
 \[
 L_\la (x,p)=L\big( J_\la (x,p),p\big) +\frac{\| x-J_\la (x,p)\|^2}{2\la}+\frac{\la\| p\|^2}{2}
\]
So (\ref{sol1}) becomes
\begin{equation}
\int_0^T  L\big(v_\la (t),\Lambda x_{\lambda}(t)+\dot x_\la (t)\big)
   +\frac{\| x_\la (t)-v_\la (t)\|^2}{2\la}
   +\frac{\la\|\Lambda x_{\lambda}(t)+\dot x_\la (t)\|^2}{2}\, dt
   +\ell\big( x_\la (0),x_\la (T)\big) =0
 \end{equation}
 where $v_\la (t)=J_\la\big(x_\la (t),\Lambda x_{\lambda}(t)+\dot x_\la (t)\big)$.
 Using  (\ref{lambda.derivative}) and  (\ref{sol2}) that  for all $t$,
\begin{equation}
 \label{sol5}
-\Lambda x_{\lambda}(t)-\dot x_\la (t)\in \partial_1L_\la \big(x_\la (t),\Lambda x_{\lambda}(t)+\dot x_\la (t)\big)
=\frac{x_\la (t)-v_\la (t)}{\la}
\end{equation}
Setting $t=0$ in  (\ref{sol2}) we get
\[\big(-\Lambda x_{0}-\dot x_\la (0),0\big)\in \partial L_\la \big( 0,\Lambda x_{0}-\dot x_\la (0)\big)\]
So  we can apply Lemma 6.3.4 to get
\[
 \|\Lambda x_{0}+\dot x_\la (0)\|_H\le C(T) \quad \hbox{\rm for all $\la >0$.}\]
This combined with (\ref{sol3}) gives
\[\|\dot x_\la (t)\|\le C(x_{0},T)\quad \hbox{\rm for all $\la >0$.}\]
It follows that $(x_\lambda)_\lambda$ is bounded in $A_H^2$ 
  and therefore there is a subsequence $(x_{\lambda_{j}})_{j}$ such that 
$x_{\lambda_{j}} (\cdot )\rightharpoonup x (\cdot )$ weakly  in $A_H^2$ and hence strongly in $C([0,T], H)$. Since $\Lambda$ is completely continuous, we also have that $\Lambda x_{\lambda_{j}}$ converges to  $\Lambda x $ in $L^2_{X^*}$
  It follows that $(\Lambda x_\lambda -\dot x_\la)$ is bounded in $L^2_{H}$ and therefore again from 
 (\ref{sol5})  we obtain
\begin{eqnarray}
\frac{\|x_\la -v_\la \|_{L^2_H}^2}{\la}\to 0
\end{eqnarray}
 It follows that 
\begin{eqnarray}
v_\la (\cdot )\rightharpoonup 
x(\cdot )\mbox{ in }L_H^2
\end{eqnarray}
and since clearly
\begin{eqnarray}
\la \frac{\|\dot x_\la (t)\|^2}{2}\to 0\mbox{ uniformly},
\end{eqnarray}
the above combine to yield that as $\la\to 0$ in  (\ref{sol1}) we get
\begin{eqnarray*}
I_{L,\ell, \Lambda}\big(x \big) =
   \int_0^T L\big(x(t), \Lambda x(t)+\dot x(t)\big)\, 
   dt+\ell\big(  x (0), x(T)\big)\le 0.
\end{eqnarray*}
Since  $I_{L,\ell, \Lambda}(x)\ge 0$ for all $x \in A_H^2$,  it follows that 
$
0=I_{L,\ell}(x)
   =\inf_{x \in A_H^2} I_{L,\ell}.
$

 \begin{proposition}
\label{Asym flow}
Let $X\subset H \subset X^*$ be an evolution triple and let $\phi : X\to \R$ be a  convex, lower semi-continuous such that for some $C>0$, $m,n>1$ we have the following growth condition:
\begin{equation}
\label{growth}
\frac{1}{C}\left( {\| x\|}_X^m-1\right)\le\phi (x)\le C\left( {\| x\|}_{X}^n+1\right).
\end{equation}
 Let $A:D(A)\subset X\to X^*$ be a skew-adjoint operator modulo boundary operators $(b_1, b_2): D(b_1,b_2)\subset X \to H\times H$ and let $\Lambda: D(\Lambda)\subset X \to X^*$ be a regular conservative operator such that $D(B)\cap D(b_1, b_2)\subset D(\Lambda)$. 
 For  $T >0$ and $v_0\in D(B)\cap D(b_1, b_2)$, define the following functional on  $L^2_X([0,T])$ by
 {\small 
 \begin{eqnarray*}
 I(u) &=& \int_0^T\left\{\phi(x(t))+\phi^*(-\Lambda x(t)-A x(t) - \dot{x}(t))+  \frac{1}{2} (|b_{1}x(t)|^{2} + |b_{2}x(t)|^{2})-2\langle b_1v_0, b_1x(t)\rangle +|b_1v_0|^{2})\right\} dt\\
 && + \frac{1}{2}(|x(0)|^{2}+|x(T)|^{2})-2\langle x(0), v_{0} \rangle +|v_{0}|^{2}, 
\end{eqnarray*}
}
 whenever $x \in S=\{ {\cal X}_2; x(t)\in D(B)\cap D(b_1, b_2)\subset D(\Lambda)\}$ and $+\infty$ elsewhere. 
Then there exists a path $v \in S $ such that
\begin{equation}
 I(x)=\inf\limits_{x\in L^2_X([0,T])}I(x)=0.
\end{equation}
 \begin{eqnarray}
 \label{eqn:1000}
-\dot{v}(t) - A v(t) -\Lambda  v(t) &\in&\partial\phi( v(t))\quad \mbox{ for a.e. }t\in [0,T] \\
 v(0)& = &v_0. \nonumber\\
b_1( v(t))&=& b_1(v_0).
\end{eqnarray}
 \end{proposition}
\noindent{\bf Proof:} The Lagrangian 
\[L(t, x,p) :=
\begin{cases}
{\phi(x) +\phi^*(-Ax - p) +   \frac{1}{2} (|b_{1}x|^{2} + |b_{2}x|^{2})-2\langle b_1v_0, b_1x\rangle +|b_1v_0|^{2}\quad {\rm if} \ x\in D(A)\cap D(b_1, b_2)}\cr {{+\infty} \quad \quad \quad \quad \quad \quad \quad \quad \quad \quad \quad \quad \quad \quad \quad \quad \quad \quad \quad \quad  {\rm elsewhere}}\cr
\end{cases} \]
is an ASD Lagrangian on $X\times X^*$ by Propostion 2.5.  . Let  $\ell$ be 
\[
\ell (r,s)=\frac{1}{2}(|r|^{2}+|s|^{2})-2\langle r, v_{0} \rangle +|v_{0}|^{2}, 
\]
and lift $\Lambda$ to a regular conservative operator $\tilde \Lambda$ from $L^2_X([0,T])$ into $L^2_{X^*}([0,T])$ by setting $(\tilde \Lambda x)(t)=\Lambda (x(t))$ which is defined on $\tilde D=\{x\in L^2_X([0,T]); x(t) \in D(\Lambda) \, \, {\rm a.e.}\}$. It is easy to check that all the conditions of Theorem 6.1 are satisfied by $L$, $\ell$, $B$ and $\tilde \Lambda$, hence  there exists $ v\in {\cal X}_2 $ such that $I(v)=0$. We obtain
\begin{eqnarray*}
0 &=& \int_0^T\big(\phi(v(t))+\phi^*(-\Lambda v(t)-A v(t) - \dot{v}(t))+ \langle v(t), \Lambda v(t)+A v(t) + \dot{v}(t)\rangle\big)\, dt\\
&&\quad + \frac{1}{2}\int_0^T\Vert b_1(v(t))-\gamma (t)\Vert^2_{H}\, dt\\
&&\quad +  \frac{1}{2}\Vert v(0)-v_0\Vert_H^2.
\end{eqnarray*}
The result follows from the fact that each term above is non-negative.

\end{document}